\pdfoutput=1
\documentclass[toc]{mathprint} %
\usepackage{rutarmacros}
\usepackage{macros}
\title[Two-scale branching functions]{Two-scale branching functions and inhomogeneous attractors}

\author
  {Vilma Orgoványi}
  {Department of Stochastics, Institute of Mathematics, Budapest University of Technology and Economics, Műegyetem rkp.\ 3., H-1111 Budapest, Hungary}
  {orgovanyi.vilma@gmail.com}

\author
  {Alex Rutar}
  {
      Department of Mathematics and Statistics,
      University of Jyväskylä,
      P.O.\ Box 35 (MaD),
      FI-40014 University of Jyväskylä,
      Finland
  }
  {alex@rutar.org}

\begin{document}
\begin{abstract}
    We introduce the notion of a \emph{two-scale branching function} associated with an arbitrary metric space, which encodes the lower and upper box dimensions as well as the Assouad spectrum.
    If the metric space is quasi-doubling, this function is approximately Lipschitz.
    We fully classify the attainable Lipschitz two-scale branching functions, which gives a new proof of the classification of Assouad spectra due to the second author.

    We then study inhomogeneous self-conformal sets satisfying standard separation conditions.
    We show that the two-scale branching function of the attractor is given explicitly in terms of the two-scale branching function of the condensation set and the Hausdorff dimension of the homogeneous attractor.
    In particular, this gives formulas for the lower box dimension and the Assouad spectrum of the attractor.
\end{abstract}

\section{Introduction}

Of all of the notions of fractal dimension, the lower and upper box dimensions are perhaps the easiest to formulate.
They are defined by
\begin{align*}
    \dimlB E &= \liminf_{r\to 0} \frac{\log \sup_{x\in E}N_r(B(x,1))}{\log(1/r)}\\
    \dimuB E &= \limsup_{r\to 0} \frac{\log \sup_{x\in E}N_r(B(x,1))}{\log(1/r)}.
\end{align*}
Here, $E$ is an arbitrary non-empty metric space, $B(x,r)$ denotes the open ball with radius $r$, and for $F\subset E$, $N_r(F)$ is the least number of open balls of radius $r$ required to cover $F$.

Beyond the lower box dimension, there are a number of fractal dimensions which encode finer scaling properties of the metric space $E$.
In this paper, we initiate the study of a \emph{two-scale branching function} associated with the metric space $E$, defined as follows.
Fix $\Delta = \{(u,v): 0\leq v \leq u\}$, and we equip $\Delta$ with the $L^1$ metric $d((u,v),(u',v')) = |u-u'| + |v-v'|$.
We then define for $(u,v)\in\Delta$
\begin{equation*}
    \beta_E(u,v) \coloneqq \log \sup_{x\in E} N_{2^{-u}}(B(x, 2^{-v})).
\end{equation*}

Here, and throughout the paper, the logarithm is in base $2$.
The function $\beta_E$, which we call two-scale branching function corresponding to the metric space $E$, takes values in the extended positive reals $[0, \infty]$.

Clearly $\beta_E$ encodes the box dimensions; it also tracks how the size of the set $E$ varies as a function of the scale.
This additional information can be very important.
For instance, it plays a key role in the lower box dimensions infinitely generated iterated function systems \cite{arxiv:2406.12821}, where the information provided by the lower and upper box dimensions alone is insufficient.
An analogue of the function $u\mapsto\beta_E(u, 0)$ has also played a key role in almost every recent paper in multi-scale incidence geometry; for a very incomplete selection, see \cite{zbl:1428.28005,zbmath:8038252,arxiv:2308.08819,arxiv:2301.10199} and the many references therein.
The key property here is that (for subsets of Euclidean space) the function $\beta_E(\cdot, 0)$ is an increasing Lipschitz function, and many problems can be reduced to an asymptotic study of Lipschitz functions.

The function $\beta_E$ also encodes the \emph{Assouad spectrum}, introduced in \cite{zbl:1390.28019}.
The Assouad spectrum is defined for $0 \leq \theta < 1$ by
\begin{equation*}
    \dimAs\theta E = \limsup_{u \to\infty}\frac{\beta_E(u, \theta u)}{u(1-\theta)}.
\end{equation*}
The Assouad spectrum is known to satisfy various functional inequalities; a full classification of the functions which can appear as the Assouad spectrum of $E$ (for subsets of Euclidean space) was given in \cite{zbmath:07937992}.

The main contributions of this paper are twofold.
\begin{itemize}
    \item First, we study the properties of the two-scale branching function $\beta_E$, with weak assumptions on $E$ (namely, that $E$ is \emph{quasi-doubling}, as defined in the next section).
        We will see that the function $\beta_E$ is fully characterized by a few simple properties, which allows one to reduce certain dimension-theoretic problems to problems which only concern the two-scale branching function.
        This handles a wide variety of technical issues in a unified way, and generalizes known results concerning Assouad spectra to a wider class of metric spaces.
        See for example \cref{ic:swap}.

    \item Secondly, to give a concrete example of the function $\beta_E$ in practice, we give a comprehensive study of the two-scale branching function $\beta_E$ of the attractor of an \emph{inhomogeneous iterated function system}.
        For background, we refer to the reader to the recent survey \cite{arxiv:2407.08404}.
        A brief introduction to inhomogeneous attractors is as follows (for more detail, see \cref{s:inhomog-2}).
        If $X$ is a compact metric space, $F\subset X$ (called the \emph{condensation set}) is compact, and $\mathcal{F}=\{f_i:i\in\mathcal{I}\}$ is a finite collection of strictly contracting Lipschitz maps from $X$ to $X$, the \emph{inhomogeneous attractor} is the unique non-empty compact set $\Lambda\subset X$ satisfying the invariance relation
        \begin{equation*}
            \Lambda = F \cup \bigcup_{i\in\mathcal{I}}f_i(\Lambda).
        \end{equation*}
        Under the additional assumptions that the maps $f_i$ are (essentially) conformal and not too concentrated in space, we will give an explicit formula for the function $\beta_\Lambda$ in terms of the function $\beta_F$ and the similarity dimension of the IFS.
        In particular, this gives precise formulas for the lower box dimension of $\Lambda$ and the Assouad spectrum of $\Lambda$.
\end{itemize}
For the remainder of the introduction, we introduce these two topics in more detail and state our main results.

\subsection{Lipschitz approximation of two-scale branching functions}
Our first key result is that the two-scale branching function can be approximated by $\alpha$-Lipschitz functions, where $\alpha$ is the quasi-Assouad dimension of the metric space.

We begin by describe appropriate spaces of Lipschitz two-scale branching functions.
\begin{definition}\label{d:B-def}
    Let $\mathcal{B}$ denote the set of functions $\psi \colon \Delta \to [0, \infty]$ with the following properties: for all $0\leq v\leq w\leq u$,
    \begin{enumerate}[nl,r]
        \item\label{e:triv} $\psi(u,u) = 0$,
        \item\label{e:subadd} $\psi(u, v) \leq \psi(u, w) + \psi(w, v)$,
        \item\label{e:incr} $\psi$ is increasing in the first variable: $\psi(u, v)\geq \psi(w, v)$, and
        \item \label{e:decr} $\psi$ is decreasing in the second variable: $\psi(u, v)\geq \psi(u, w)$.
    \end{enumerate}
    Then for $\alpha \geq 0$, we let $\mathcal{B}(\alpha)\subset\mathcal{B}$ denote those functions which are $\alpha$-Lipschitz.
    For simplicity, we also write $\mathcal{B} = \mathcal{B}(\infty)$.
\end{definition}
Now, let $E$ be a non-empty metric space with two-scale branching function $\beta_E$.
The initial observation is that $\beta_E\in\mathcal{B}$.
Properties \cref{e:triv}, \cref{e:incr}, and \cref{e:decr} are almost immediate from the definition.
The subadditivity property \cref{e:subadd} follows by fixing a cover for $B(x, 2^{-v})$ with balls $B(y, 2^{-w})$, and then covering each ball of radius $2^{-w}$ at scale $2^{-u}$:
\begin{equation*}
    N_{2^{-u}}(B(x, 2^{-v})) \leq \sup_{y\in E}N_{2^{-u}}(B(y,2^{-w}))\cdot N_{2^{-w}}(B(x, 2^{-v})).
\end{equation*}

Handling the $\alpha$-Lipschitz condition is more subtle.
We will see that the equivalent condition on the set $E$ can be phrased in terms of the \emph{quasi-Assouad dimension}, which is defined as follows:
\begin{align*}
    \dimqA E = \lim_{\theta \nearrow 1} \inf\Bigl\{s \geq 0&:\exists C > 0\, \forall 0<r \leq r^\theta \leq R < 1\\*
                                                     &\sup_{x\in E}N_r(B(x,R) \cap E) \leq C \left(\frac{R}{r}\right)^s\Bigr\}.
\end{align*}
We say that $E$ is \emph{quasi-doubling} if $\dimqA E < \infty$.
If the metric space is \emph{doubling} (for instance, if $E\subset \R^d$ for some $d \in \N$), then it is automatically quasi-doubling.

Unpacking definitions, if $\dimqA E = \alpha$, then $\beta_E$ is bounded in an approximate sense: for all $\varepsilon>0$, there exists a $C_\varepsilon>0$ such that for all $0\leq \theta \leq 1-\varepsilon$ and $u\geq 0$,
\begin{equation*}
    \beta(u, \theta u) \leq C_\varepsilon + (u - \theta u)(\alpha + \varepsilon).
\end{equation*}

Because of the error terms, it is unreasonable to expect to expect that $\beta_E$ will be literally a member of $\mathcal{B}(\alpha)$.
On the other hand, this is indeed the case up to an error term which is invisible to quantities such as the Assouad spectrum.
\begin{itheorem}\label{it:Ba}
    Let $E$ be a non-empty metric space.
    Then
    \begin{equation*}
        \dimqA E = \min\{\alpha \geq 0:\exists \psi\in\mathcal{B}(\alpha)\text{ s.t.\ } \beta_E(u,v) = \psi(u,v) + o(u)\}.
    \end{equation*}
    Conversely, if $d\in\N$, $\alpha\in [0,d]$, and $\psi\in\mathcal{B}(\alpha)$, then there exists a compact set $E\subset\R^d$ such that
    \begin{equation*}
        \beta_E(u,v) = \psi(u,v) + O(\log(u-v+1)).
    \end{equation*}
\end{itheorem}
In other words, up to error $o(u)$, $\mathcal{B}(\alpha)$ is precisely the class of two-scale branching functions associated to metric spaces with quasi-Assouad dimension at most $\alpha$.
The reader is encouraged to regard the space $\mathcal{B}(\alpha)$ as the space of two-scale branching functions in $\mathcal{B}$ which can be approximated by $\alpha$-Lipschitz functions.

Let us illustrate the technical utility of \cref{it:Ba} with a short proof.
We first recall the definition of the upper Assouad spectrum: for $0 \leq \theta < 1$, it is given by
\begin{equation*}
    \dimuAs\theta E = \limsup_{u \to\infty}\sup_{0\leq \lambda \leq \theta}\frac{\beta(u, \lambda u)}{u(1-\lambda)}.
\end{equation*}
We obtain the following generalization of \cite[Theorem~2.1]{zbl:1410.28008}.
\begin{icorollary}\label{ic:swap}
    Let $E$ be a non-empty metric space with $\dimqA E < \infty$.
    Then for $0 \leq \theta < 1$,
    \begin{equation}\label{e:swap}
        \dimuAs \theta E = \sup_{0\leq\lambda\leq\theta}\dimAs\lambda E.
    \end{equation}
\end{icorollary}
\begin{proof}
    Applying \cref{it:Ba}, get a function $\psi\in\mathcal{B}(\alpha)$ such that $\beta_E(u,v) = \psi(u,v) + o(u)$.
    Next, by definition of the upper Assouad spectrum, get an increasing sequence $0 < u_n \to \infty$ and $\lambda_n \in [0, \theta]$ such that $\lambda_n \to \lambda \in [0, \theta]$ and
    \begin{equation*}
        \dimuAs\theta E = \lim_{n\to\infty} \frac{\beta_E(u_n, \lambda_n u_n)}{u_n(1-\lambda_n)} = \lim_{n\to\infty} \frac{\psi(u_n, \lambda_n u_n)}{u_n(1-\lambda_n)}.
    \end{equation*}
    Since $\psi$ is $\alpha$-Lipschitz, for $n\in\N$,
    \begin{equation*}
        |\psi(u_n, \lambda_n u_n) - \psi(u_n, \lambda)| \leq \alpha u_n|\lambda_n - \lambda|.
    \end{equation*}
    Since $\theta < 1$, $(1-\lambda_n)^{-1} \to (1-\lambda)^{-1}$, so
    \begin{equation*}
        \lim_{n\to\infty} \left\lvert\frac{\psi(u_n, \lambda_n u_n)}{u_n(1-\lambda_n)} - \frac{\psi(u_n, \lambda)}{u_n(1-\lambda)}\right\rvert = 0.
    \end{equation*}
    Therefore $\dimAs\lambda E \geq \dimuAs\theta E$ as required.
\end{proof}
Essentially, \cref{e:swap} is a statement about swapping a limit and a supremum, which is justified since the family of functions $\mathcal{B}(\alpha)$ is uniformly equicontinuous.

\subsection{Normalized limits and the Assouad spectrum}\label{ss:norm-intro}
Next, we turn our attention to a certain normalized limit associated with a two-scale branching function $\beta\in\mathcal{B}$.
This limit encodes the Assouad spectrum of the corresponding set and is unchanged by $o(u)$ error terms.

We begin with a space of limits of two-scale branching functions.
\begin{definition}\label{d:G-def}
    We let $\mathcal{G}$ denote the set of functions $\gamma\colon [0,1] \to [0, \infty]$ such that:
    \begin{enumerate}[nl,r]
        \item $\gamma(1) = 0$,
        \item $\theta\mapsto\gamma(\theta)$ is decreasing, and
        \item\label{i:G-subadd} For all $\lambda, \theta\in[0,1]$,
            \begin{equation*}
                \gamma(\lambda\theta) \leq \gamma(\theta) + \theta\gamma(\lambda).
            \end{equation*}
    \end{enumerate}
    For $\alpha \geq 0$, we let $\mathcal{G}(\alpha)\subset\mathcal{G}$ denote the $\alpha$-Lipschitz elements of $\mathcal{G}$.
\end{definition}
The space $\mathcal{G}(\alpha)$ arises naturally from the following normalized limit.
If $\psi\in\mathcal{B}(\alpha)$ and $u > 0$, the function
\begin{equation*}
    \psi_u(\theta)\coloneqq \frac{\psi(u, \theta u)}{u}
\end{equation*}
is a decreasing $\alpha$-Lipschitz function from $[0,1]$ to $[0, \alpha]$ with $\psi_u(1) = 0$.
In particular, it makes sense to define, in the sense of uniform convergence,
\begin{equation*}
    \Gamma(\psi) = \limsup_{u\to\infty} \psi_u.
\end{equation*}
Re-writing \cref{d:B-def}~\cref{e:subadd} using this notation,
\begin{equation}\label{e:rescale-def}
    \psi_u(\lambda \theta) - \psi_u(\theta) \leq \theta \psi_{\theta u}(\lambda).
\end{equation}
Taking a limit supremum in $u$ yields \cref{d:G-def}~\cref{i:G-subadd}.
The other properties are easy to check, so one can verify that $\Gamma$ maps functions in $\mathcal{B}(\alpha)$ to $\mathcal{G}(\alpha)$.

The function $\Gamma$ encodes the Assouad spectrum in terms of the two-scale branching function.
If $E$ has $\dimqA E = \alpha < \infty$ and $\beta_E(u,v) = \psi(u,v)+o(u)$ for $\psi\in\mathcal{B}(\alpha)$, then
\begin{equation*}
    (1-\theta)\dimAs\theta E = \Gamma(\psi)(\theta).
\end{equation*}
The key property of $\Gamma$ is the following.
\begin{itheorem}\label{it:surj}
    For $\alpha \geq 0$, the map $\Gamma\colon\mathcal{B}(\alpha) \to \mathcal{G}(\alpha)$ is a continuous order-preserving surjection.
\end{itheorem}
This is the symbolic analogue of the classification of Assouad spectra from \cite{zbmath:07937992} and, when combined with the classification part of \cref{it:Ba}, gives a new proof of the classification.

Finally, we note that the two-scale branching function is preserved by the quasi-Lipschitz maps introduced in \cite{zbl:1145.28007}.
This generalizes a similar observation for the quasi-Assouad dimension from \cite{zbl:1345.28019}.
See \cref{s:ql-inv} for more detail.

\subsection{Monotone subspaces and inhomogeneous attractors}
In the final two sections, we focus on the two-scale branching functions associated with inhomogeneous attractors.
Given an inhomogeneous attractor $\Lambda_F$ with condensation set $F$, we will given an explicit formula for $\beta_\Lambda$ (up to $o(u)$ error) in terms of $\beta_F$ and a certain critical exponent.
We will assume that the IFS is \emph{minimally distorting} (see \cref{d:md}) and \emph{asymptotically bounded} (see \cref{d:asymp}).

We defer discussion to \cref{s:inhom-assump}, but for the moment, we note that these assumptions are satisfied with critical exponent $\dimH\Lambda_\varnothing = \dimB\Lambda_\varnothing$ in the following cases:
\begin{enumerate}[nl]
    \item $X \subset \R$, the maps $f_i$ are similarity maps satisfying the \emph{exponential separation condition}, and $\dimH \Lambda_\varnothing < 1$; this follows from the results in \cite{zbl:1426.11079}.
    \item $X$ is a compact connected subset of $\R^d$, the maps $f_i$ are conformal maps with bounded distortion, and the IFS satisfies the open set condition; see \cite{zbl:0852.28005}.
        For definitions and a proof of the minimal distortion property, see \cref{s:conf}.
\end{enumerate}

The explicit formula is given by the function $\Phi_h\colon\mathcal{B}(\alpha)\to\mathcal{B}(\alpha)$ defined for $0\leq h\leq \alpha$ by
\begin{equation*}
    \Phi_h(\psi)(u,v) = \max_{0 \leq z \leq u}
    \begin{cases}
        \psi(u-z,v-z) &: 0 \leq z \leq v,\\
        h(z-v) + \beta(u-z, 0) &: v \leq z \leq u.
    \end{cases}
\end{equation*}
First, in \cref{s:inhomog} we study the abstract properties of the function $\Phi_h$ and the space $\Phi_h(\mathcal{B}(\alpha))$.
We will see that $\Phi_h$ is a projection onto the space $\Bh(\alpha)$ defined as follows:
\begin{definition}\label{d:Bsub}
    For $0 \leq h \leq \alpha$, we let $\Bh(\alpha)\subset\mathcal{B}(\alpha)$ denote the functions $\psi\in\mathcal{B}(\alpha)$ which satisfy the following additional properties:
    \begin{enumerate}[nl,r]
        \item\label{i:diag-mono} For all $(u,v)\in\Delta$ and $z\geq 0$, $\psi(u+z,v+z)\geq\psi(u,v)$.
        \item For all $(u,v)\in\Delta$, $\psi(u,0)-\psi(v,0) \geq h(u-v)$.
    \end{enumerate}
\end{definition}
Then, in \cref{s:inhomog-2}, we prove that $\beta_{\Lambda} = \Phi_h(\psi)$ (up to $o(u)$ error), where $\psi$ is any $\alpha$-Lipschitz approximation of $\beta_F$ with error at most $o(u)$.

This is summarized in the following theorem.
\begin{itheorem}\label{it:D}
    Let $X$ be a metric space with $\dimqA X = \alpha < \infty$.
    Suppose $F$ has approximate two-scale branching function $\psi \in \mathcal{B}(\alpha)$.
    Then $\Phi_h(\psi) \in \Bh(\alpha)$ and moreover
    \begin{equation*}
        \Phi_h(\beta)=\inf\{\xi\in \Bh(\alpha): \xi\geq \beta\}.
    \end{equation*}
    In particular, $\Phi_h\colon \mathcal{B}(\alpha)\to \overline{\mathcal{B}}_h(\alpha)$ is surjective, idempotent, and the identity on its range.

    Moreover, suppose $X$ and $F$ are compact and $\Lambda\subset X$ is an inhomogeneous attractor with condensation set $F$, and that the IFS is minimally distorting and asymptotically bounded.
    Then $\Lambda$ has approximate two-scale branching function $\Phi_h(\psi)$.
\end{itheorem}
The first part is stated and proved as an independent theorem in \cref{t:mono-sub-B}, and the proof of the second part can be found in \cref{t:in-dim}.

This explicit formula also gives a complete classification of the attainable two-scale branching functions of inhomogeneous attractors.

To conclude, we show how this result gives a concrete formula for the Assouad spectrum of $\Lambda_F$ in terms of the upper Assouad spectrum of $F$ and the critical exponent $h$.
In particular, we will prove that the map $\Phi_h$ commutes with the limiting operator $\Gamma$ defined in the previous section.

For notational convenience, for $\theta\in[0,1]$, for $\gamma\in\mathcal{G}(\alpha)$ we define $\Psi(\gamma)\in\mathcal{G}(\alpha)$ by the rule
\begin{equation}
    \Psi(\gamma)(\theta) \coloneqq
    \begin{cases}
        \frac{\gamma(\theta)}{1-\theta} &: 0 \leq \theta < 1,\\
        \lim_{\theta\nearrow 1}\frac{\gamma(\theta)}{1-\theta} &: \theta = 1.
    \end{cases}
\end{equation}
That the limit at $1$ exists is well-known from \cite[Proposition~2.3]{zbmath:07937992}.

Next, we define the function $\Omega_h(\gamma)$, which is analogous to $\Phi_h$ but instead defined on the space of limits $\mathcal{G}(\alpha)$.
\begin{definition}
    For $\gamma \in \mathcal{G}(\alpha)$ and $\theta\in[0, 1]$, we write
    \begin{equation*}
        \Omega_h(\gamma)(\theta) \coloneqq (1-\theta)\cdot\max\left\{h, \max_{0\leq \lambda \leq \theta}\Psi(\lambda)\right\}.
    \end{equation*}
\end{definition}
Recall that $\Phi_h$ is a projection onto the subspace $\Bh(\alpha)\subset\mathcal{B}(\alpha)$: the analogous subspace for $\Omega_h$ is the space $\Gh(\alpha)\subset\mathcal{G}(\alpha)$ defined as follows.
\begin{definition}\label{d:Gsub}
    For $0 \leq h\leq \alpha$, let $\Gh(\alpha)\subset \mathcal{G}(\alpha)$ denote the set of functions $\gamma\in \mathcal{G}(\alpha)$  which satisfies:
    \begin{enumerate}[nl,r]
        \item\label{i:incr} $\theta\mapsto\gamma(\theta)/(1-\theta)$ is increasing, and
        \item $\gamma(0) \geq h$.
    \end{enumerate}
\end{definition}
In fact, the subadditivity property is implied by \cref{i:incr} in the above definition; we give the details in \cref{s:alt}.

The above can be summarized in the following theorem.
\begin{itheorem}\label{it:E}
    Fix notation as in \cref{it:D}.
    Then the maps $\Gamma|_{\Bh(\alpha)}\colon \Bh(\alpha) \to \Gh(\alpha)$ and $\Omega_h\colon\mathcal{G}(\alpha)\to \Gh(\alpha)$ are surjective and following diagram commutes:
    \begin{equation*}
        \begin{tikzcd}
            \mathcal{B}(\alpha) \arrow{r}{\Phi_h} \arrow[swap]{d}{\Gamma} & \Bh(\alpha) \arrow{d}{\Gamma} \\%
            \mathcal{G}(\alpha) \arrow{r}{\Omega_h}& \Gh(\alpha)
        \end{tikzcd}
    \end{equation*}

    In particular, the inhomogeneous attractor $\Lambda$ has Assouad spectrum
    \begin{equation*}
        \dimAs\theta \Lambda = \Psi\circ\Gamma\circ\Phi_h(\beta)(\theta)=\Psi\circ\Omega_h\circ\Gamma(\beta)(\theta) = \max\{h, \dimuAs\theta F\}.
    \end{equation*}
\end{itheorem}
For surjectivity of $\Omega_h$, see \cref{t:mono-sub}.
Then for surjectivity of $\Gamma$ restricted to the subspace $\Gh(\alpha)$, and commutativity of the diagram, see \cref{t:mono-sub-equiv}.

\subsection{Lower box dimension}
Finally, we briefly note an application to lower box dimension, and a relationship with the main result of \cite{arxiv:2406.12821}.
This generalizes the main results of \cite{zbl:1342.28013}, where partial results were obtained in terms of a certain \emph{covering regularity exponent}.
Related results (but in somewhat different contexts) can be found in \cite{zbl:1423.28017,zbl:1403.37024}.

Let us begin by noting the following application of \cref{it:D}.
Given a quasi-doubling metric space $E$ with approximate two-scale branching function $\psi\in\mathcal{B}(\alpha)$, let $g_E(u) = \psi(u,0)$ for $u\geq 0$.

If $E$ is bounded, then
\begin{equation*}
    g_E(u) = \log N_{2^{-u}}(E) + o(u).
\end{equation*}
In other words, the function $g_E$ encodes the upper and lower box dimension of $E$.
We call the function $g_E$ the \emph{average branching function} associated with $E$.
\begin{icorollary}\label{ic:lower-box}
    Let $X$ be a compact metric space with $\dimqA X = \alpha < \infty$, let $F\subset X$ be compact, and let $\Lambda\subset X$ be an inhomogeneous attractor with condensation set $F$ which is minimally distorting and asymptotically bounded.
    Let $F$ have average branching function $g_F$.
    Then $\Lambda$ has average branching function
    \begin{equation*}
        g_\Lambda(u) = \sup_{0\leq z \leq u}\left(g_F(z) + h(u-z)\right) + o(u).
    \end{equation*}
\end{icorollary}
This formula is the same formula as stated in \cite[Theorem~C]{arxiv:2406.12821}, where $\Lambda$ is the limit set of an infinitely generated self-conformal IFS, and $F$ is the set of fixed points of the maps.

Now, we fix some notation from \cite[Theorem~D]{arxiv:2406.12821}: given $0\leq s \leq t \leq\alpha$ and $0\leq h\leq\alpha$, write
\begin{equation*}
    \mathcal{D}(h,s,t,\alpha)\coloneqq \begin{cases}
        \{h\} &: t\leq h,\\
        \left[\max\{h,s\}, h+\frac{(t-h)(\alpha-h)s}{\alpha t-h s}\right] &: t>h.
    \end{cases}
\end{equation*}
The following corollary follows by identical arguments as given in \cite{arxiv:2406.12821}.
\begin{icorollary}
    Fix notation as in \cref{ic:lower-box}.
    Then the following hold:
    \begin{enumerate}[nl,r]
        \item We have $\dimlB\Lambda=\dimuB\Lambda$ if and only if
            \begin{equation*}
                \dimuB F \leq \max\{h, \dimlB F\}.
            \end{equation*}
        \item We have $\dimlB \Lambda\in\mathcal{D}(h, \dimlB F, \dimuB F, \alpha)$.
        \item Suppose $0\leq h\leq \alpha$, $0\leq s \leq t \leq \alpha$, and $\zeta\in\mathcal{D}(h,s,t,\alpha)$.
            Then for any $d\in\N$ with $d\geq \alpha$, there is an inhomogeneous self-similar set satisfying the open set condition on $\R^d$ with limit set $\Lambda$ and condensation set $F$ such that $\dimH\Lambda = h$, $\dimlB F = s$, $\dimuB F = t$, and $\dimlB \Lambda = \zeta$.
    \end{enumerate}
\end{icorollary}
Essentially, infinitely generated self-conformal sets and inhomogeneous self-conformal sets are structurally quite similar, at least on average.
On the other hand, the formulas for the Assouad spectrum established in \cref{it:E} \emph{do not} hold for infinitely self-conformal sets.
Such examples can be found in \cite{zbl:1558.28002}.

\subsection{Further thoughts: \texorpdfstring{$L^p$}{Lp} generalizations}
Finally, we note a potential generalization of the two-scale branching function, and note a few properties.

Let $E$ be a metric space, let $1 \leq p\leq\infty$ and define for $(u,v)\in\Delta$
\begin{equation*}
    \beta_E^p(u,v) = \log \norm{N_{2^{-u}}\bigl(B(x, 2^{-v})\bigr)}_p.
\end{equation*}
We are intentionally imprecise with how to define $\norm{\cdot}_p$, but let us give a sketch of one sensible choice.
Fix a non-empty finite subset $C\subset E$ with the property that $B(x, 2^{-v})\cap B(y, 2^{-v})=\varnothing$ for $x\neq y\in C$.
Then we can equip $C$ with normalized counting measure, and take the $L^p$-norm relative to this measure, over points $x\in C$.
Then we write
\begin{equation*}
    \norm{N_{2^{-u}}\bigl(B(x, 2^{-v})\bigr)}_p
\end{equation*}
to denote the supremum of this quantity over all finite subsets $C$ of $E$ with the separation property.

In the case $p=\infty$ this is exactly the two-scale branching function.
On the other hand, when $p=1$,
\begin{equation*}
    \beta_E^1(u,v) = \beta_E(u, 0) - \beta_E(v,0) + O(1).
\end{equation*}
If $E$ is bounded, this is just the average covering number between two scales:
\begin{equation*}
    \beta_E^1(u,v) = \log\left(\frac{N_{2^{-u}}(E)}{N_{2^{-v}}(E)}\right) + O(1).
\end{equation*}

For uniformly branching sets (see \cref{r:unif}), $\beta_E^1 = \beta_E^\infty$.
In particular, it is reasonable to interpret the function $\beta_E^p$ for $1 < p < \infty$ as a measurement of uniformity.

The natural question is therefore if there is an interesting theory for $1<p<\infty$, and if such a theory would have applications.
We leave such considerations to further work for the interested reader.

\section{A two-scale branching function}
\subsection{Branching functions}
Recall from the introduction that
\begin{equation*}
    \Delta = \{(u,v): 0\leq v \leq u\},
\end{equation*}
which we equip with the $L^1$ metric: $d((u,v),(u',v')) = |u-u'| + |v-v'|$.
In \cref{d:B-def} we defined the space $\mathcal{B}$ of two-scale branching functions with subspace $\mathcal{B}(\alpha)\subset \mathcal{B}$ consisting of those functions which are $\alpha$-Lipschitz.

We begin with a few equivalent ways to guarantee the $\alpha$-Lipschitz property.
\begin{lemma}\label{l:lip-char}
    Let $\psi\in \mathcal{B}$ and $\alpha \geq 0$.
    Then the following are equivalent.
    \begin{enumerate}[nl,a]
        \item\label{i:lip} $\psi$ is $\alpha$-Lipschitz.
        \item\label{i:lip-1} For all $v \geq 0$, the function $u\mapsto \psi(u,v)$ is $\alpha$-Lipschitz.
        \item\label{i:lip-2} For all $u \geq 0$, the function $v\mapsto \psi(u,v)$ is $\alpha$-Lipschitz.
        \item\label{i:bd} For all $(u,v) \in \Delta$,
            \begin{equation*}
                \psi(u, v) \leq \alpha(u-v).
            \end{equation*}

    \end{enumerate}
\end{lemma}
\begin{proof}
    It is immediate that \cref{i:lip}$\Rightarrow$\cref{i:lip-1} and \cref{i:lip}$\Rightarrow$\cref{i:lip-2}, that  \cref{i:lip-1} and \cref{i:lip-2} together imply \cref{i:lip}.
    Moreover, \cref{i:lip-1}$\Rightarrow$\cref{i:bd} and \cref{i:lip-2}$\Rightarrow$\cref{i:bd} since $\psi(u,u) = 0$ for all $u \geq 0$.

    Therefore, it remains to show that \cref{i:bd}$\Rightarrow$\cref{i:lip-1} and \cref{i:bd}$\Rightarrow$\cref{i:lip-2}.
    Indeed, to see \cref{i:bd}$\Rightarrow$\cref{i:lip-2}, if $0\leq v<w \leq u$ are arbitrary, then by property \cref{e:decr}, \cref{e:subadd} of \cref{d:B-def}, and assumption \cref{i:bd}
    \begin{equation*}
        0 \leq \psi(u,v) - \psi(u,w)\leq \psi(w,v)\leq \alpha(w-v).
    \end{equation*}
   The proof of \cref{i:bd}$\Rightarrow$\cref{i:lip-1} is analogous.
\end{proof}
\begin{remark}
    One can further prove that the following condition is equivalent with the above: \textit{\begin{enumerate}[label=(e)]
        \item\label{i:bd-lim} For all $u > 0$,
            \begin{equation*}
                \limsup_{v \to u}\frac{\psi(u, v)}{u - v} \leq \alpha.
            \end{equation*}
            \end{enumerate}}
   \noindent That \cref{i:bd}$\Rightarrow$\cref{i:bd-lim} is immediate and \cref{i:bd-lim}$\Rightarrow$\cref{i:lip-2} and \cref{i:bd-lim}$\Rightarrow$\cref{i:lip-1} can be proven by a mean value theorem argument (see, for example, \cite[Corollary 2.3]{zbl:1509.28005}).
\end{remark}
\subsection{Approximation by Lipschitz branching functions}
We now show that the space $\mathcal{B}(\alpha)$ provides an appropriate representation of the two-scale branching function associated with a metric space $E$.
Recall, for a metric space $E$ that the two-scale branching function is defined for $(u,v)\in\Delta$ by
\begin{equation*}
    \beta_E(u,v) = \log \sup_{x\in E}N_{2^{-u}}\bigl(B(x, 2^{-v})\bigr).
\end{equation*}

First recall the following lemma, the proof of which is straightforward and sketched in the introduction.
\begin{lemma}\label{l:B-member}
    Let $E$ have two-scale branching function $\beta$.
    Then $\beta\in\mathcal{B}$.
\end{lemma}

The remainder of this section is dedicated to the proof of the first half of \cref{it:Ba}: if $\alpha = \dimqA E < \infty$, then $\beta_E \in \mathcal{B}(\alpha)$ with an appropriate error term.

The idea of the proof is as follows.
The upper bound $\beta_{E}(u,v)\leq \dimqA(E)(u-v)+o(u)$ can be readily proved from the definition of $\dimqA E$ using \cref{l:lip-char}.
The lower bound is somewhat more difficult: we must construct an appropriate function $\xi\in\mathcal{B}(\alpha)$.
The main idea is to use a certain minimal extension of $\alpha$-Lipschitz functions.
Given $b\geq 0$, we first construct a one variable $\alpha$-Lipschitz function which approximates $\beta_E$ along the line segment ${(u, b)}_{u\geq b}$ and then extend this to a function in $\mathcal{B}(\alpha)$ using \cref{l:approx-to-B}.
The minimality property of this extension ensures that the resulting function is bounded above by $\beta_E$.
Then taking the supremum of these extensions for all $b\geq 0$ gives the desired element of $\mathcal{B}(\alpha)$.

We now proceed with the implementation of this proof.
We begin by showing that the function $\beta_E$ is approximately $\alpha$-Lipschitz.
\begin{lemma}\label{l:approx-lip}
    Let $E$ be a non-empty metric space with $\alpha = \dimqA E < \infty$.
    Then for all $(u,v)\in\Delta$,
    \begin{equation*}
        \beta_E(u,v) \leq \alpha(u-v) + o(u).
    \end{equation*}
\end{lemma}
\begin{proof}
    It suffices to show for all $\varepsilon > 0$ there is a constant $C_\varepsilon \geq 1$ such that for all $0<r \leq R < 1$,
    \begin{equation}\label{e:approx-sub}
        \sup_{x\in E}N_r(B(x,R)) \leq C_\varepsilon r^{-\varepsilon} \left(\frac{R}{r}\right)^\alpha.
    \end{equation}
    Let $\theta < 1$ be sufficiently large so that $(1-\theta)(\alpha + \varepsilon) \leq \varepsilon$.
    We consider two cases depending on the value of $R$.
    If $R \geq r^\theta$, then by definition of the quasi-Assouad dimension there is a constant $C_\varepsilon$ (depending on $E$, $\theta$, and $\varepsilon$) so that
    \begin{equation*}
        N_r(B(x,R)) \leq C_\varepsilon \left(\frac{R}{r}\right)^{\alpha + \varepsilon} \leq C_\varepsilon r^{-\varepsilon} \left(\frac{R}{r}\right)^\alpha.
    \end{equation*}
    In the second inequality we just use that $R < 1$.
    Otherwise, if $R \leq r^\theta$, since $B(x,R) \subset B(x, r^\theta)$,
    \begin{equation*}
        N_r(B(x,R)) \leq N_r(B(x, r^\theta)) \leq C_\varepsilon\left(\frac{r^\theta}{r}\right)^{\alpha + \varepsilon} \leq C_\varepsilon r^{-(1-\theta)(\alpha + \varepsilon)} \leq C_\varepsilon r^{-\varepsilon}\left(\frac{R}{r}\right)^\alpha
    \end{equation*}
    where the last line follows since $(R/r)^\alpha \geq 1$.
    Since $x\in E$ was arbitrary, the claim in \cref{e:approx-sub} follows.

    Taking logarithms and substituting the definition of $\beta_E$, the desired claim follows.
\end{proof}
We now show that the approximate $\alpha$-Lipschitz inequality in \cref{l:approx-lip} implies good approximation by a function in $\mathcal{B}(\alpha)$.

In order to do this, we introduce a particular family of functions in $\mathcal{B}(\alpha)$.
We begin with a family of Lipschitz functions.
\begin{definition}\label{d:lip-0}
    Given $\alpha \geq 0$, we let $\mathcal{C}(\alpha)$ denote the set of functions $g\colon [0,\infty)\to [0,\infty)$ which are increasing, $\alpha$-Lipschitz, and have $g(0) = 0$.
\end{definition}
We now describe a process which allows us to extend a Lipschitz function to a two-parameter function $\psi \in \mathcal{B}(\alpha)$ which is minimal in the following sense.
\begin{lemma}\label{l:h-ext}
    Let $\alpha \geq 0$ and $g \in \mathcal{C}(\alpha)$.
    Let $\xi\colon \Delta \to [0,\infty)$ be given by $\xi(u,v) = g(u) - g(v)$.
    Then $\xi \in \mathcal{B}(\alpha)$.

    Moreover, if for some $b \geq 0$, $g(b) = 0$, and $\beta \in \mathcal{B}(\alpha)$ is any function with
    \begin{equation*}
        \beta(u, b) - \beta(v,b) \geq g(u) - g(v)\quad\text{for all}\quad b \leq v \leq u,
    \end{equation*}
    then $\beta \geq \xi$.
\end{lemma}
\begin{proof}
    Clearly $\xi(u,u) = 0$ for all $u\geq 0$, and $\xi(u,v)$ is increasing in $u$ and decreasing in $v$.
    Moreover, \cref{e:subadd} holds (with equality) since
    \begin{equation}\label{e:subadd-exact}
        \xi(u,v) = g(u) - g(w) + g(w) - g(v) = \xi(u,w) + \xi(w,v).
    \end{equation}
    Therefore, $\xi \in \mathcal{B}$.
    Moreover, since $g$ is $\alpha$-Lipschitz, $u \mapsto \xi(u,v)$ is $\alpha$-Lipschitz, and therefore $\xi\in\mathcal{B}(\alpha)$ by \cref{l:lip-char}.

    Finally, let $b \geq 0$, and suppose $g(b) = 0$ and $\beta \in \mathcal{B}(\alpha)$ has $\beta(a,b) = \xi(a,b)$ for all $a \geq b$.
    Suppose $(u,v)\in\Delta$.
    If $v \geq b$, then
    \begin{equation*}
        \beta(u,v) \geq \beta(u, b) - \beta(v, b) \geq g(u) - g(v) = \xi(u,v);
    \end{equation*}
    and if $v \leq b$, then since $g(b) = g(v) = 0$,
    \begin{equation*}
        \beta(u,v) \geq \beta(u, b) = \beta(u,b) - \beta(b,b) \geq  g(u) - g(b) = \xi(u,v)
    \end{equation*}
    as required.
\end{proof}
\begin{remark}\label{r:unif}
    For a fixed function $g \in \mathcal{C}(\alpha)$, the function $\xi$ is essentially the two-scale branching function associated with a \emph{uniformly branching set} (a set which can be represented by a uniformly branching tree) with branching function $g$.

    In the case $b=0$, the latter half of the lemma states that if $E$ is an arbitrary set with branching function $g_E$ and two-scale branching function $\beta_E$, then $\beta_E \geq \xi$.
    In other words, the uniformly branching sets minimize the two-scale branching functions among all sets with a given branching function.
\end{remark}
Next, we note that $\mathcal{B}(\alpha)$ and $\mathcal{C}(\alpha)$ are supremum-closed.
\begin{lemma}\label{l:sup-closed}
    Let $\mathcal{F}\subset \mathcal{B}(\alpha)$ (resp.\ $\mathcal{F}\subset \mathcal{C}(\alpha)$) be non-empty.
    Then the pointwise supremum $f(x)=\sup{g(x): g\in \mathcal{F}}$ is in $\mathcal{B}(\alpha)$ (resp.\ $\mathcal{C}(\alpha)$).
\end{lemma}
\begin{proof}
    A direct computation shows that the maximum of two elements in $\mathcal{B}(\alpha)$ (resp.\ $\mathcal{C}(\alpha)$) is still an element of $\mathcal{B}(\alpha)$.
    Then the desired claim follows by a compactness argument using, for instance, the Arzelà–Ascoli theorem.
\end{proof}
Now, we show that an approximate $\alpha$-Lipschitz property is sufficient for good approximation by a function in $\mathcal{B}(\alpha)$.
We state this as a separate result with an explicit error term since this may be useful in applications where better approximation is required.
\begin{proposition}\label{l:approx-to-B}
    Let $\eta \colon [0,\infty) \to [0,\infty)$ be increasing, let $\alpha \geq 0$, and suppose $\beta \in \mathcal{B}$ has
    \begin{equation}\label{e:apl}
        \beta(u,v) \leq \alpha (u-v) + \eta(u)
    \end{equation}
    for all $(u,v) \in \Delta$.
    Then there exists a function $\psi\in\mathcal{B}(\alpha)$ such that
    \begin{equation*}
         \beta(u,v)- \eta(u) \leq \psi(u,v) \leq \beta(u,v)+ \eta(u).
    \end{equation*}
\end{proposition}
\begin{proof}
    We will construct the function $\psi \in \mathcal{B}(\alpha)$ as a supremum of functions $\xi_{b} \in \mathcal{B}(\alpha)$ for $b \geq 0$.
    The supremum is in $\mathcal{B}(\alpha)$ by \cref{l:sup-closed}.

    Fix $b \geq 0$ and consider the function
    \begin{equation*}
        g_b = \sup\{f \in \mathcal{C}(\alpha): f(a) \leq \beta(a,b)\text{ for all }a\geq b\}.
    \end{equation*}
    Since $\beta \geq 0$, the supremum is non-empty, so by \cref{l:sup-closed}, $g_b \in \mathcal{C}(\alpha)$, and $g_b(a) \leq \beta(a,b)$ for all $a\geq b$ by construction.
    We next show that
    \begin{equation}\label{e:gb}
        g_b(a) \geq \beta(a,b) - \eta(a)\qquad\text{for all}\qquad a \geq b.
    \end{equation}
    Let $a \geq b$ be fixed.
    If $\beta(a,b) \leq \eta(a)$, there is nothing to prove.
    Otherwise, consider the function $f$ such that $f(u) = \beta(a,b) - \eta(a)$ for all $u \geq a$, and for $u \leq a$,
    \begin{equation*}
        f(u) = \max\{f(a) - \alpha(a - u), 0\}.
    \end{equation*} 

    We will prove that $f \in \mathcal{C}(\alpha)$ and $f(u) \leq \beta(u,b)$ for all $u \geq b$, from which it follows that $g_b(a)\geq f(a)$ and hence \cref{e:gb} holds.
    First, $f \geq 0$ is certainly increasing and $\alpha$-Lipschitz and $f(0)=0$, since $f(a)-\alpha(a-u)=0$ at $u=a-f(a)/\alpha\geq 0$ by \cref{e:apl}.
    Moreover, if $b \leq u \leq a$, is such that $f(u) > 0$, then since $\beta \in \mathcal{B}$,
    \begin{align*}
        f(u) - \beta(u,b) &= \beta(a,b) - \beta(u,b) - \eta(a) - \alpha(a-u)\\
                          & \leq \beta(a, u) - \eta(a) - \alpha(a-u)\\
                          & \leq \alpha(a-u) + \eta(a) - \eta(a) - \alpha(a-u)\\
                          & = 0
    \end{align*}
    where we used \cref{e:subadd} (and that $a\geq u$) and \cref{e:apl}.
    Finally, if $u \geq a$, since $f(a) \leq \beta(a,b)$, $\beta$ is increasing in the first coordinate, and $f$ is constant on $[a,\infty)$, we conclude that $f(u) \leq \beta(u,b)$ for all $u \geq b$.
    Since $a \geq b$ was arbitrary, we have shown that \cref{e:gb} holds.

    Finally we extend $g_b$; for $b\geq 0$, let $\xi_b(u,v) = g_b(u) - g_b(v)$.
    Of course, $\xi_b \in \mathcal{B}(\alpha)$ by \cref{l:h-ext}.
    
    Now we prove that $\xi_b(u, v)\leq \beta(u, v)+\eta(u)$ for all $(u, v)\in \Delta$.
    A preliminary observation is from \cref{e:gb}, recalling that $g_b(b) = 0$, that
    \begin{equation}\label{e:beta-xi}
        \beta(u,b) \geq \xi_b(u, b) = g_b(u) \geq \beta(u,b) - \eta(u)\qquad\text{for all}\qquad u \geq b.
    \end{equation}
    Hence, for an arbitrary $(u,v) \in \Delta$ we have for $v \geq b$, using \cref{e:subadd}, \cref{e:beta-xi}, and \cref{e:apl},
    \begin{equation*}
        \beta(u,v) \geq \beta(u,b) - \beta(v, b)  \geq \xi_b(u,b) - \xi_b(v,b) -  \eta(u) = \xi_b(u,v) \eta(u).
    \end{equation*}
    Otherwise, if $v \leq b$, since $\beta\in\mathcal{B}$,
    \begin{equation*}
        \beta(u,v) \geq \beta(u, b) \geq \xi_b(u, b) = \xi_b(u,v)
    \end{equation*}
    since $g_b(v) = g_b(b) = 0$.

    Conversely, for $(u, v)\in \Delta$, the function $\xi_v$ satisfies $\xi_v(u,v) \geq \beta(u,v) - \eta(u)$ by \cref{e:gb}.
    Therefore, the $\psi = \sup\{\xi_{b}: b \geq 0\}$ satisfies the required properties.
\end{proof}
We are now in position to establish the main result of this section.
\begin{theorem}\label{t:approx-Lip}
    Let $E$ be a non-empty metric space with two-scale branching function $\beta_E$.
    Then there exists $\psi \in \mathcal{B}(\dimqA E)$ such that for all $(u,v) \in \Delta$,
    \begin{equation*}
        \beta_E(u,v) = \psi(u,v) + o(u).
    \end{equation*}
    In fact,
    \begin{equation*}
        \dimqA E = \min\{\alpha \geq 0:\exists \psi\in\mathcal{B}(\alpha)\text{ s.t.\ } \beta_E(u,v) = \psi(u,v) + o(u)\}.
    \end{equation*}
\end{theorem}
\begin{proof}
    Write $\alpha = \dimqA E$.
    If $\alpha = \infty$, then $\beta_E \in \mathcal{B} = \mathcal{B}(\infty)$ and there is nothing to prove.
    Otherwise, if $\alpha < \infty$, combining \cref{l:approx-lip} and \cref{l:approx-to-B}, it follows that there is a $\psi\in\mathcal{B}(\alpha)$ such that for all $(u,v)\in\Delta$, $\beta_E(u,v) = \psi(u,v) + o(u)$.

    For the second part of the theorem, suppose $\alpha \geq 0$ is such that $\beta_E(u,v) = \psi(u,v) + o(u)$ where $\psi \in \mathcal{B}(\alpha)$.
    To complete the proof, it suffices to show that $\dimqA E \leq \alpha$.

    Let $\eta\colon[0,\infty) \to [0,\infty)$ be an increasing function with $\lim_{u\to\infty} u^{-1}\eta(u) = 0$ such that $|\beta_E(u,v) - \psi(u,v)| \leq \eta(u)$.

    If $\alpha = \infty$, there is nothing to prove.
    Otherwise, let $\varepsilon>0$ be arbitrary: we prove that $\dimqA E \leq \alpha + \varepsilon$.
    Let $0<\theta < 1$ be arbitrary.
    Rephrasing the definition of $\dimqA E$ in terms of $\beta_E$, it suffices to prove that there is a constant $A \geq 0$ so that for all $(u,v)\in\Delta$ with $v \leq \theta u$ that
    \begin{equation}\label{e:beta}
        \beta_E(u,v) \leq A + (u-v)(\alpha + \varepsilon).
    \end{equation}

    Since $u^{-1}\eta(u)$ converges to $0$, there is a $u_0 \geq 0$ so that for all $u \geq u_0$, $\eta(u) \leq u(1-\theta)\varepsilon$.
    Therefore since $\psi\in\mathcal{B}(\alpha)$, for all $(u,v) \in \Delta$ with $u\geq u_0$ and $v \leq \theta u$,
    \begin{align*}
        \beta_E(u,v) &\leq \psi(u,v) + \eta(u)\\
                   &\leq \alpha(u-v) + \eta(u)\\
                   &\leq \alpha(u-v) + u(1-\theta) \varepsilon\\
                   & \leq (u-v)(\alpha + \varepsilon).
    \end{align*}
    In the above computation, we only used $v \leq \theta u$ in the last line.
    On the other hand, for $u \leq u_0$, since $\beta_E$ is increasing in $u$ and decreasing in $v$, the first two lines of the above computation shows that
    \begin{equation*}
        \beta_E(u,v) \leq \beta_E(u_0, 0) \leq \alpha u_0 + u_0(1-\theta)\varepsilon < \infty.
    \end{equation*}
    Therefore taking $A = u_0 \alpha + u(1-\theta)\varepsilon$, we conclude that \cref{e:beta} holds for all $(u,v) \in \Delta$ with $v \leq \theta u$ and.
\end{proof}

\subsection{Attainable branching functions}\label{ss:attain}
In this section, we will prove a converse to \cref{t:approx-Lip}: namely, that every $\psi\in\mathcal{B}(\alpha)$ appears as the approximate two-scale branching function of a compact subset of Euclidean space.
We begin with a lemma to approximate Lipschitz functions by integer-valued step functions.
First, some notation: for $\alpha\geq 0$, let $\mathcal{S}(\alpha)$ denote the set of functions $\eta\colon\N_0\to\alpha\N_0$ for which $\eta(0) = 0$ and $\eta(n) - \eta(n-1)\in\{0,\alpha\}$ for all $n\in\N$.

The space $\mathcal{S}(\alpha)$ is a discretized version of $\mathcal{C}(\alpha)$, which we recall from \cref{d:lip-0} denotes the increasing $\alpha$-Lipschitz functions satisfying $g(0)=0$.
\begin{lemma}
    Let $\alpha \geq 0$ and $g\in\mathcal{C}(\alpha)$.
    Then there exists an $\eta\in\mathcal{S}(\alpha)$ such that for all $n\in\N$, $g(n) - \alpha  < \eta(n) \leq g(n)$.
\end{lemma}
\begin{proof}
    We construct $\eta$ inductively as follows.
    Let $\eta(0)=0$, and inductively set
    \begin{equation*}
        \eta(n+1) = \begin{cases}
            \eta(n)+\alpha &: g(n+1) - \eta(n) \geq \alpha\\
            \eta(n) &: g(n+1) - \eta(n) < \alpha.
        \end{cases}
    \end{equation*}
    Certainly $\eta\in\mathcal{S}(\alpha)$.

    We now prove by induction that $g(n) - \alpha < \eta(n) \leq g(n)$ for all $n\in\N_0$.
    Clearly $\eta(0) = g(0)$; suppose the claim holds for $n\in\N_0$.

    If $g(n+1) - \eta(n) \geq \alpha$, then $g(n+1)-\eta(n+1)=g(n+1)-\eta(n)-\alpha\geq 0$ by the definition of $\eta$, and
    \begin{equation*}
        g(n+1)-\eta(n)-\alpha\leq g(n)+\alpha-\eta(n)-\alpha
    \end{equation*}
    since $g$ is $\alpha$-Lipschitz.
    Otherwise, if $g(n+1) - \eta(n) < \alpha$ then $g(n+1)-\eta(n+1)=g(n+1)-\eta(n)<\alpha$ by the definition of $\eta$, and $g(n+1)-\eta(n)\geq g(n)-\eta(n)\geq 0$.
    In either case, we have $g(n+1) - \alpha < \eta(n+1) \leq g(n+1)$ as desired.
\end{proof}
Our next lemma is analogous to \cite[Lemma~3.7]{zbl:1509.28005}.
\begin{lemma}\label{l:homog}
    Let $0\leq \alpha \leq d$ and $g\in\mathcal{C}(\alpha)$.
    Suppose $m\in\N_0$ is such that $g(m) = 0$.
    Then there exists a compact set $E\subset [0,2^{-m}]^d$ such that
    \begin{equation*}
        \beta_E(u,v) = g(u) - g(v) + O(1)
    \end{equation*}
    with implicit constants depending only on $d$.
\end{lemma}
\begin{proof}
    First, get $\eta\in\mathcal{S}(d)$ such that for all $n\in\N$, $g(n) - d < \eta(n) \leq g(n)$.
    We now use the sequence $\eta$ to define a set by subdivision as follows.
    Let $\mathcal{D}_n$ denote the set of closed grid-aligned dyadic cubes of side-length $2^{-n}$ contained in $[0,1]^d$.
    We inductively define subsets of $[0,1]^d$ as unions of dyadic cubes in $\mathcal{D}_n$.
    Begin with the unit cube $[0,1]^d$.
    Now, suppose we have constructed some family of cubes $\mathcal{Q}_n\subset\mathcal{D}_n$ for $n\in\N_0$.
    If $\eta(n)=0$, replace each $Q\in\mathcal{Q}_n$ by a new cube $Q'\in\mathcal{Q}_{n+1}$ sharing the bottom-left corner; otherwise, replace $Q$ by the set of all cubes $Q'\in\mathcal{D}_{n+1}$ for which $Q'\subset Q$.
    Finally, we define
    \begin{equation*}
        E = \bigcap_{n=0}^\infty\bigcup_{Q\in\mathcal{Q}_n} Q.
    \end{equation*}
    Clearly $E$ is non-empty and compact, and since $g(m) = 0$, $E\subset [0,2^{-m}]^d$.

    Finally, let $0 \leq v \leq u$ be arbitrary: we estimate $\beta_E(u,v)$.
    Let $k,n\in\N_0$ be such that $v-1 < k \leq v \leq u \leq n < u+1$.
    By construction, if $x\in E$ is arbitrary, then $B(x, 2^{-v})\cap E$ intersects $O(1)$ dyadic cubes of side-length $2^{-k}$ in $\mathcal{Q}_k$, and each cube $Q\in\mathcal{Q}_k$ contains $2^{\eta(n) - \eta(k)}$ dyadic cubes in $\mathcal{Q}_n$.
    Therefore
    \begin{equation*}
        \log N_{2^{-u}}(B(x, 2^{-v})\cap E) \leq \eta(n) - \eta(k) + O(1) = g(u) - g(v) + O(1).
    \end{equation*}
    It is clear that this bound is sharp up to $O(1)$ error since the dyadic cubes in $\mathcal{Q}_n$ are disjoint, as required.
\end{proof}
We now obtain the main result of this section, and with it complete the proof of \cref{it:Ba}.
The fundamental building block of this proof is similar to the idea in the proof of \cref{l:approx-to-B}: for each fixed $b\geq 0$, we construct a set from a uniformly branching dyadic tree with branching numbers defined by the function $\psi$ on the line $\{(u,b)\}_{u\geq b}$, and then take a union over $b\in\N_0$ to obtain the desired set.
\begin{theorem}\label{t:attain}
    Let $d\in\N$ and $0\leq \alpha \leq d$.
    Let $\psi\in\mathcal{B}(\alpha)$ be arbitrary.
    Then there exists a compact set $E\subset\R^d$ such that $\beta_E(u,v) = \psi(u,v) + O(\log(u-v+1))$.
\end{theorem}
\begin{proof}
    Let us begin by approximating the function $\psi(u,v)$ as a maximum of functions provided by \cref{l:h-ext}.

    For each $b\in\N_0$, define
    \begin{equation*}
        g^{b}(u) =
        \begin{cases}
            0 &: u\leq b,\\
            \psi(u, b) &: u\geq b.
        \end{cases}
    \end{equation*}
    Then with $\xi^b(u,v) = g^b(u) - g^b(v)$, by \cref{l:h-ext}, $\xi^b \leq \psi$ and $\xi^b(u,b) = \psi(u,b)$ for all $u\geq b$.
    In particular, $\xi\coloneqq \max\{\xi^b:b\in\N_0\}$ satisfies
    \begin{equation*}
        \psi(u,v) = \xi(u,v) + O(1)
    \end{equation*}
    with implicit constant depending only on the Lipschitz constant of $\psi$.

    Now, for each $b\in\N_0$, apply \cref{l:homog} to obtain a set $E_b\subset [0, 2^{-b}]^d$ such that
    \begin{equation*}
        \beta_{E_b}(u,v) = g^b(u) - g^b(v) + O(1) = \xi^b(u,v) + O(1).
    \end{equation*}
    Finally, write
    \begin{equation*}
        E = \{(0,\ldots,0)\} \cup \bigcup_{b\in \N_0}(E_b + 2^{-b+2}).
    \end{equation*}
    By definition, $E$ is a closed subset of $[0,5]^d$.
    Moreover, $\dist(E_b, E\setminus E_b) > 2^{-b}$.

    Since $E_b\subset E$ for all $b\in\N_0$, $\beta_E \geq \beta_{E_b}$ and
    \begin{equation*}
        \beta_E(u,v) \geq \xi(u,v) + O(1) = \psi(u,v) + O(1).
    \end{equation*}
    Conversely, let $0 \leq v \leq u$ be arbitrary and let $x\in E$.
    Let $n\geq u$ be minimal such that $n\in\N$.
    Then $E_b\subset[0,2^{-n}]^d$ for all $b\geq n$ and therefore
    \begin{equation*}
        N_{2^{-u}}\left(\bigcup_{b=n}^\infty (E_b + 2^{-b+2})\right) = O(1).
    \end{equation*}
    Since the remaining sets $E_b + 2^{-b+2}$ are $2^{-b}$-separated, the ball $B(x, 2^{-v})$ intersects at most $u-v+1$ of the sets, and therefore
    \begin{align*}
        \beta_E(u,v) &\leq \log \left((u-v+1)\cdot 2^{\max\{\beta_{E_b}(u,v): 0\leq b \leq n\}}\right) + O(1)\\
                     &\leq \max\{\xi^b(u,v) : 0\leq b\leq n\} + O(\log(u-v+1))\\
                     &\leq \psi(u,v) + O(\log(u-v+1))
    \end{align*}
    as claimed.
\end{proof}

\subsection{Normalized limits of the branching function}\label{ss:norm-lims}
In this section, we describe a certain normalized limit associated with a two-scale branching function; we will see that this limit encodes the Assouad spectrum.
We also prove \cref{it:surj}.

Recall for $\alpha \geq 0$ the definition of $\mathcal{G}(\alpha)$ and $\psi_u$ for $\psi\in\mathcal{B}(\alpha)$ and $u>0$ from \cref{ss:norm-intro}, as well as the function $\Gamma\colon\mathcal{B}(\alpha)\to\mathcal{G}(\alpha)$.

Let us first note that the convergence is uniform in the following sense.
\begin{lemma}\label{l:gamma-unif}
    Let $\alpha \geq 0$ and let $\beta$ be a function such that $\beta(u,v) = \psi(u,v) + o(u)$ for $\psi\in\mathcal{B}(\alpha)$.
    Let $\gamma = \Gamma(\beta)$.
    Then $\gamma\in\mathcal{G}(\alpha)$, and
    \begin{equation*}
        \limsup_{u\to\infty}\sup_{0\leq \theta \leq 1}\left(\frac{\beta(u, \theta u)}{u} - \gamma(\theta)\right) \leq 0.
    \end{equation*}
\end{lemma}
\begin{proof}
    It suffices to show that
    \begin{equation*}
        \limsup_{u\to\infty}\sup_{0\leq \theta \leq 1}\left(\gamma(\theta) - \psi_u(\theta)\right) \leq 0.
    \end{equation*}
    First, for $u>0$, define a new function
    \begin{equation*}
        g_u = \sup_{v\geq u}\psi_v.
    \end{equation*}
    One can check that $g_u$ is $\alpha$-Lipschitz and takes values in $[0,\alpha]$.
    Moreover, by definition,
    \begin{equation*}
        \gamma(\theta) = \lim_{u\to\infty}g_u(\theta)
    \end{equation*}
    pointwise.
    Since $g_u$ is a monotonic sequence, by the Arzelà--Ascoli theorem, the convergence is uniform in $\theta$.
    But $\psi_u(\theta)\leq g_u(\theta)$ for all $u>0$ and $\theta\in[0,1]$, so the claim follows.
\end{proof}
It is proven in \cite{zbmath:07937992} that $\mathcal{G}(\alpha)$ is in bijection with the set of functions which appear as the Assouad spectrum of a subset of Euclidean space.
In particular, with \cref{t:attain} in mind, it follows implicitly that $\Gamma$ is a surjective function.
We now give a substantially simpler explicit proof.
In fact, we establish the existence of certain maximal inverse for $\Gamma$.
\begin{proposition}\label{p:Gamma-inv}
    Let $\alpha \geq 0$ and let $\gamma\in\mathcal{G}(\alpha)$ be arbitrary.
    Define $\psi(u,v) = u\cdot\gamma(v/u)$.
    Then:
    \begin{enumerate}[nl,r]
        \item $\psi\in\mathcal{B}(\alpha)$,
        \item $\Gamma(\psi) = \gamma$, and
        \item If $\beta\in\mathcal{B}(\alpha)$ has $\Gamma(\beta) \leq \gamma$, then $\beta(u,v) \leq \psi(u,v) + o(u)$.
    \end{enumerate}
    In particular, $\Gamma\colon\mathcal{B}(\alpha)\to\mathcal{G}(\alpha)$ is surjective.
\end{proposition}
\begin{proof}
    Let us first check that $\psi\in\mathcal{B}$.
    Clearly $\psi(u,u) = 0$ since $\gamma(1) = 0$.
    Moreover, $\psi$ is increasing in $u$ and decreasing in $v$ since $\gamma$ is decreasing.
    Then to check subadditivity, if $0\leq v\leq w \leq u$ are arbitrary, by \cref{d:G-def}~\cref{i:G-subadd},
    \begin{equation*}
        \psi(u,v) = u\cdot\gamma(v/u) \leq u\left(\gamma(w/u) + (w/u)\gamma(v/w)\right) = \psi(u,w) + \psi(w,v)
    \end{equation*}
    as required.

    Next, since $\gamma$ is $\alpha$-Lipschitz, it follows that $v\mapsto\psi(u,v)$ is $\alpha$-Lipschitz, so $\psi$ is $\alpha$-Lipschitz by \cref{l:lip-char}.
    Therefore $\psi\in\mathcal{B}(\alpha)$.

    Finally, $\psi$ is defined so that for all $u>0$ and $\theta\in[0,1]$, $\psi_u(\theta) = \gamma(\theta)$.
    Therefore it is clear that $\Gamma(\psi)=\gamma$ and $\psi \geq \beta$ for all $\beta$ for which $\Gamma(\beta) = \gamma$.
\end{proof}
With this, we complete the proof of \cref{it:surj} since continuity and the order-preserving property are immediate from the definition of $\Gamma$.

\subsection{Assouad spectrum and quasi-Assouad dimension}
Now, let $E$ be quasi-doubling with $\alpha = \dimqA E$.
We introduce another space of functions, which is in bijection with the space $\mathcal{G}(\alpha)$.
\begin{definition}\label{d:A-def}
    For $\alpha\geq 0$, we set
    \begin{equation*}
        \mathcal{A}(\alpha) = \left\{\varphi \in C([0,1]):\theta\mapsto(1-\theta)\varphi(\theta) \in \mathcal{G}(\alpha)\right\}.
    \end{equation*}
\end{definition}
It is proven in \cite{zbmath:07937992} that $\mathcal{A}(\alpha)$ is precisely the set of functions attainable as Assouad spectra of subsets of Euclidean space with quasi-Assouad dimension at most $\alpha$.

Moreover, it is proven in \cite[Proposition~2.3]{zbmath:07937992} that
\begin{equation}\label{e:gamma-lim-1}
    \lim_{\theta \nearrow 1} \frac{\gamma(\theta)}{1-\theta}= \sup_{\theta \in [0,1)}\frac{\gamma(\theta)}{1-\theta}.
\end{equation}
With \cref{e:gamma-lim-1} in mind, recall the map $\Psi\colon \mathcal{G}(\alpha) \to \mathcal{A}(\alpha)$ defined for $0\leq \theta < 1$
\begin{equation*}
    \Psi(\gamma)(\theta) = \frac{\gamma(\theta)}{1-\theta}
\end{equation*}
and extended to $\theta = 1$ by continuity is homeomorphism.

The function $\Gamma$ computes the Assouad spectrum of $E$: for $0\leq \theta < 1$, if $\beta_E$ is the two-scale branching function of $E$, then
\begin{equation*}
    \dimAs\theta E = \Psi\circ\Gamma(\beta_E)(\theta).
\end{equation*}
Since $\gamma$ is $\alpha$-Lipschitz, for all $0 \leq \theta \leq 1$, $\gamma(\theta) \leq \alpha(1-\theta)$.
Moreover, the smallest number $\alpha\geq 0$ for which such an inequality can hold is exactly the quasi-Assouad dimension.
This result was established in \cite{zbl:1410.28008} under the additional assumption that $E$ is a doubling metric space, and we give the short proof for quasi-doubling spaces.
\begin{proposition}
    Let $E$ be a non-empty metric space with two-scale branching function $\beta_E$, and suppose $\alpha = \dimqA E < \infty$.
    Then $\alpha = \Psi\circ\Gamma(\beta_E)(1)$.
\end{proposition}
\begin{proof}
    By \cref{e:gamma-lim-1}, we may set
    \begin{equation*}
        \alpha_0 = \lim_{\theta \nearrow 1}\frac{\gamma(\theta)}{1-\theta} = \sup_{0\leq \theta < 1}\frac{\gamma(\theta)}{1-\theta}
    \end{equation*}
    We already saw above that $\alpha_0 \leq \alpha$.
    For the converse inequality, we will show that for all $(u,v) \in \Delta$ 
    \begin{equation*}
        \beta_E(u,v) \leq \alpha_0(u-v) + o(u).
    \end{equation*}
    From this, the result follows by combining \cref{l:approx-to-B} and the second part of \cref{t:approx-Lip}.

    To prove the above claim, recall from \cref{l:gamma-unif} that
    \begin{equation*}
        \limsup_{u\to\infty}\sup_{0 \leq \theta \leq 1}\left(\frac{\beta_E(u,\theta u)}{u} - \gamma(\theta)\right) \leq 0.
    \end{equation*}
    Multiplying by $u$ and writing $v = \theta u$,
    \begin{equation*}
        \beta_E(u, v) \leq u\cdot\gamma(\theta) + o(u) \leq u\cdot\alpha_0(1-\theta) + o(u) = \alpha(u-v) + o(u),
    \end{equation*}
    as required.
\end{proof}

\subsection{Quasi-Lipschitz invariance of branching functions}\label{s:ql-inv}
To conclude this section, we show that the two-scale branching function is preserved by a certain family of maps which generalize bi-Lipschitz maps between metric space.

Since we are willing to accept error terms which grows sufficiently slowly, we introduce a bit more notation.
\begin{definition}\label{d:err}
    We say that a function $\eta$ is an \emph{error function} if $\eta\colon[0,\infty)\to[0,\infty)$ is increasing and satisfies $\lim_{u\to\infty}u^{-1}\eta(u) = 0$.
    We denote the space of such functions by $\mathcal{E}$.
    We then say that a metric space $E$ is \emph{$\alpha$-dimensional with error $\eta$} if $\eta\in\mathcal{E}$ and there is a function $\psi\in\mathcal{B}(\alpha)$ such that for all $(u,v)\in\Delta$,
    \begin{equation*}
        |\psi(u,v) - \beta_E(u,v)| \leq \eta(u).
    \end{equation*}
\end{definition}
If $E$ is quasi-doubling, then there is some finite $\alpha$ and function $\eta\in\mathcal{E}$ such that $E$ is $\alpha$-dimensional with error $\eta$; this is \cref{t:approx-Lip}.
\begin{definition}
    Let $X$ and $Y$ be metric spaces and let $\eta\in\mathcal{E}$.
    We say that a surjective map $f\colon X\to Y$ is \emph{$\eta$-uniform} if for all $0\leq v\leq u$,
    \begin{equation*}
        |\beta_Y(u,v) - \beta_X(u,v)| \leq \eta(u).
    \end{equation*}
\end{definition}
It turns out that the quasi-Lipschitz maps introduced in \cite{zbl:1145.28007} provide a general class of uniform maps.
Let us first recall the definition:
\begin{definition}
    Let $X$ and $Y$ be metric spaces and let $\eta\in\mathcal{E}$.
    We say that a surjective map $f\colon X\to Y$ is \emph{$\eta$-quasi-Lipschitz} if for all $x_1,x_2\in X$,
    \begin{equation*}
        |\log d\bigl(f(x_1), f(x_2)\bigr) - \log d(x_1,x_2)| \leq \eta\bigl(-\log d(x_1, x_2)\bigr).
    \end{equation*}
\end{definition}
If $f$ is bi-Lipschitz, then it is $\eta$-quasi-Lipschitz for a bounded function $\eta$.
On the other hand, Hölder functions are certainly not quasi-Lipschitz.

It is straightforward to verify that if $f$ is quasi-Lipschitz, then $f$ is invertible and moreover $f^{-1}$ is also quasi-Lipschitz with a (potentially different) function $\tilde\eta$ depending only on $\eta$.
We say that $X$ and $Y$ are \emph{$\eta$-quasi-Lipschitz equivalent} if $f$ and $f^{-1}$ are both $\eta$-quasi-Lipschitz.

The normalization of $\eta$ is chosen for the following reason: if $f\colon X\to Y$ is $\eta$-quasi-Lipschitz, for all $x\in X$ and $w\geq 0$,
\begin{equation}\label{e:ball-map}
    B(f(x), 2^{-w-\eta(w)}) \subseteq f\bigl(B(x, 2^{-w})\bigr) \subseteq B(f(x), 2^{-w+\eta(w)}).
\end{equation}

In \cite{zbl:1345.28019} where the notion of quasi-Assouad dimension was introduced, it was shown that quasi-Assouad dimension is a quasi-Lipschitz invariant.
We now show generally that the two-scale branching function is quasi-Lipschitz invariant.
\begin{proposition}\label{p:pres}
    Let $X$ and $Y$ be quasi-doubling metric spaces which are $\eta$-quasi-Lipschitz equivalent.
    Let $\alpha<\infty$ and $\kappa\in\mathcal{E}$ be such that $X$ and $Y$ are $\alpha$-dimensional with error $\kappa$.
    Then
    \begin{equation*}
        \beta_Y(u,v) = \beta_X(u,v) + o_{\alpha, \eta, \kappa}(u).
    \end{equation*}
\end{proposition}
\begin{proof}
    Let $v\geq 0$ be such that $v\geq \eta(v)$, let $u\geq v$, and let $x\in X$ be arbitrary.
    Then by \cref{e:ball-map}, $f(B(x, 2^{-v}))\subset B(f(x), 2^{-v + \eta(v)})$.
    Applying the definition of $\beta_Y$, we can cover $B(f(x), 2^{-v+\eta(v)})$ by $2^{\beta_Y(u+\eta(u), v-\eta(v))}$ balls of radius $2^{-u-\eta(u)}$, say $\{B(y_n, 2^{-u-\eta(u)})\}_{n=1}^m$.
    But then $f(B(f^{-1}(y_n), 2^{-u}))\supset B(y_n, 2^{-u - \eta(u)})$.
    Therefore $\{B(f^{-1}(y_n), 2^{-u})\}_{n=1}^m$ is a cover for $B(x, 2^{-v})$ so
    \begin{equation*}
        \beta_X(u,v) \leq \beta_Y(u+\eta(u), v - \eta(v)) \leq \beta_Y(u,v) + 2\alpha \eta(u) +\kappa(u) + \kappa(u+\eta(u))
    \end{equation*}
    by definition of $\kappa$.

    The same argument with $X$ and $Y$ swapped yields the other inequality, as required.
\end{proof}

\section{Branching functions of inhomogeneous self-conformal sets}\label{s:inhomog}
We now turn our attention to the two-scale branching functions of inhomogeneous self-conformal sets.

Let us first briefly recall the key spaces and functions introduced in \cref{d:B-def,d:G-def,d:A-def} from \cref{ss:norm-lims}:
\begin{equation*}
    \mathcal{B}(\alpha) \overset{\Gamma}\longrightarrow \mathcal{G}(\alpha)\overset{\Psi}\longleftrightarrow \mathcal{A}(\alpha).
\end{equation*}
Here, $\Psi$ is the map $\Psi(\gamma)(\theta) = \gamma(\theta)/(1-\theta)$ extended to take value at $1$ by continuity.

We also recall from \cref{d:Bsub} in the introduction the definition of $\Bh(\alpha)\subset\mathcal{B}(\alpha)$.
We will see that the first condition in \cref{d:Bsub} corresponds to a monotonicity property of the Assouad spectrum, and the second condition corresponds to a certain minimal $h$-dimensionality.

\subsection{Projection onto monotone subspaces}
Let $0 \leq h\leq \alpha$.
We recall the definition of the $\Phi_h$, defined for $\beta\in\mathcal{B}(\alpha)$ by the rule
\begin{equation*}
    \Phi_h(\beta)(u,v) = \max_{0 \leq z \leq u}
    \begin{cases}
        \beta(u-z,v-z) &: 0 \leq z \leq v,\\
        h(z-v) + \beta(u-z, 0) &: v \leq z \leq u.
    \end{cases}
\end{equation*}
The maximum is necessarily attained by continuity of $\beta$.
\begin{lemma}\label{l:phi-elem}
    Let $0 \leq h\leq \alpha$.
    Then $\Phi_h(\beta)\in\Bh(\alpha)$.
\end{lemma}
\begin{proof}
    It is immediate that $\xi(u,u) = 0$ for all $u \geq 0$.

    To check monotonicity, it will be convenient to introduce some notation.
    For $(u,v)\in\Delta$, let
    \begin{equation*}
        E_{u,v} = \{(u-z,v-z):0\leq z\leq u\}\cup\{(u-z, 0): v \leq z \leq u\}.
    \end{equation*}
    and define $\eta_{u,v}\colon E_{u,v} \to [0,\infty)$ by
    \begin{align*}
        &\eta_{u,v}(u-z,v-z) = \beta(u-z,v-z)\qquad\text{for}\qquad 0\leq z \leq v\\
        &\eta_{u,v}(u-z, 0) = h(z-v) + \beta(u-z, 0)\qquad\text{for}\qquad v\leq z \leq u.
    \end{align*}
    In other words, $E_{u,v}$ and $\eta_{u,v}$ are chosen so that $\Phi_h(\beta)(u,v) = \max\{\eta_{u,v}(u',v'):(u',v')\in E_{u,v}\}$.

    Now to see monotonicity along diagonals (that is, \cref{d:Bsub}~\cref{i:diag-mono}), let $(u,v)\in\Delta$ and $z\geq 0$.
    Observe that $E_{u,v}\subset E_{u+z,v+z}$ and moreover $\eta_{u+z,v+z} = \eta_{u,v}$ on $E_{u,v}$.
    Therefore $\xi(u+z,v+z) \geq \xi(u,v)$.
    Next we check that $v \mapsto \xi(u,v)$ is decreasing.
    Let $0\leq v \leq w \leq u$ be arbitrary and let $0 \leq z \leq u$.
    If $z \geq w$, then
    \begin{equation*}
        h(z-v) + \beta(u-z, 0) \geq h(z-w) + \beta(u-z,0);
    \end{equation*}
    if $w \geq z \geq v$, then since $\beta(u-z,\cdot)$ is decreasing,
    \begin{equation*}
        h(z-v) + \beta(u-z, 0) \geq \beta(u-z,w-z);
    \end{equation*}
    and if $v \geq z$, then since $\beta(u-z,\cdot)$ is decreasing,
    \begin{equation*}
        \beta(u-z,v-z) \geq \beta(u-z,w-z).
    \end{equation*}
    In particular, by definition of $\xi$, it follows that $\xi(u,v) \geq \xi(u,w)$.
    Finally, by combining the above monotonicity properties,
    \begin{equation*}
        \xi(w,v) \leq \xi(u, v+u-w) \leq \xi(u,v)
    \end{equation*}
    which is the final monotonicity property, as required.

    Next, we check subadditivity.
    Let $0\leq u\leq w\leq v$ be arbitrary, and get $0 \leq z \leq u$ which attains the maximum in the definition of $\Phi_h(\beta)$.
    First suppose $0\leq z \leq v$.
    Then since $\xi(u,w) \geq \beta(u-z,w-z)$ and $\xi(w,v) \geq \beta(w-z,v-z)$, by subadditivity of $\beta$,
    \begin{align*}
        \xi(u,v) &= \beta(u-z,v-z)\\
                 &\leq \beta(u-z,w-z) + \beta(w-z,v-z)\\
                 &\leq \xi(u,w) + \xi(w,v).
    \end{align*}
    Next suppose $v \leq z \leq w$.
    Then since $\xi(w,v) \geq \beta(w-z,0) + h(z-v)$ and $\xi(u,w) \geq \beta(u-z,w-z)$, by subadditivity of $\beta$,
    \begin{align*}
        \xi(u,v) &= h(z-v) + \beta(u-z, 0)\\
                 &\leq h(z-v) + \beta(w-z,0) + \beta(u-z,w-z)\\
                 &\leq h(z-v) + (\xi(w,v) - h(z-v)) + \xi(u,w)\\
                 &= \xi(w,v) + \xi(u,w).
    \end{align*}
    Finally, suppose $w \leq z \leq u$.
    Then since $\xi(u,w) \geq \beta(u-z, 0) + h(z-w)$ and $\xi(w,v) \geq h(w-v)$,
    \begin{align*}
        \xi(u,v) &= h(z-v) + \beta(u-z, 0)\\
                 &\leq h(z-v) + \xi(u,w) - h(z-w) + \xi(w,v) - h(w-v)\\
                 &= \xi(u,w) + \xi(w,v)
    \end{align*}
    as required.

    It remains to show that $\xi(u,0) - \xi(v,0) \geq h(u-v)$ for all $0\leq v\leq u$.
    Let $0\leq z \leq v$ be such that $\xi(v,0) = h\cdot z +\beta(v-z, 0)$.
    Then
    \begin{equation*}
        \xi(u,0) \geq h(u-v+z) + \beta(v-z, 0) = h(u-v) + \xi(v,0)
    \end{equation*}
    as required.
\end{proof}
In fact, $\Phi_h$ is a projection onto $\Bh(\alpha)$.
This is the second part of \cref{it:D}.
\begin{theorem}\label{t:mono-sub-B}
    Let $0 \leq h\leq \alpha$ and $\beta\in\mathcal{B}(\alpha)$.
    Then $\Phi_h(\beta)\in\Bh(\alpha)$, and moreover
    \begin{equation*}
        \Phi_h(\beta) = \inf\{\xi\in\Bh(\alpha): \xi \geq \beta\}.
    \end{equation*}
    In particular, $\Phi_h\colon\mathcal{B}(\alpha) \to\Bh(\alpha)$ is a surjective idempotent map which is the identity on $\Bh(\alpha)$.
\end{theorem}
\begin{proof}
    Let $\beta\in\mathcal{B}(\alpha)$ be arbitrary, and write $\xi = \Phi_h(\beta)$.

    We already saw in \cref{l:phi-elem} that $\xi\in\Bh(\alpha)$.
    Moreover, taking $z=0$ in the definition of $\Phi_h$, it is clear that $\xi \geq \beta$.
    Therefore it remains to show that if $\eta\in\Bh(\alpha)$ with $\eta \geq \beta$ is arbitrary, then $\eta \geq \xi$.
    Let $(u,v)\in\Delta$ be fixed and get $0 \leq z \leq u$ which attains the maximum in the definition of $\beta$.
    If $0 \leq z \leq v$, then using monotonicity of $\xi$ along diagonals,
    \begin{equation*}
        \xi(u,v) = \beta(u-z,v-z) \leq \eta(u-z,v-z) \leq \eta(u,v)
    \end{equation*}
    and if $v \leq z \leq u$, then using the lower $h$-Lipschitz bound of $\xi(\cdot, 0)$ followed by monotonicity along diagonals,
    \begin{equation*}
        \xi(u,v) = \beta(u-z,0) + h (z-v) \leq \xi(u-z, 0) + h(z-v) \leq \xi(u-v, 0) \leq \xi(u,v).
    \end{equation*}
    In either case, $\xi \leq \eta$, as required.
\end{proof}

\subsection{Assouad spectra associated with monotone subspaces}
In the previous section, we studied the subspace $\Bh(\alpha)$ along with the corresponding projection $\Phi_h\colon \mathcal{B}(\alpha) \to \Bh(\alpha)$.

Now, recall the limiting function $\Gamma\colon\mathcal{B}(\alpha) \to \mathcal{G}(\alpha)$, which maps the two-scale branching function to the corresponding Assouad spectrum function.
We also recall from the introduction the subspace $\Gh(\alpha)\subset \mathcal{G}(\alpha)$.

In this section, we will prove that the space $\Gh(\alpha)$ corresponds directly to $\Bh(\alpha) \subset \mathcal{B}(\alpha)$.
The correspondence is via the function $\Omega_h$ defined for $\gamma\in\mathcal{G}(\alpha)$ by the rule
\begin{equation*}
    \Omega_h(\gamma)(\theta) = (1-\theta)\cdot\max\left\{h, \max_{0\leq \theta' \leq \theta} \Psi(\gamma)(\theta')\right\}.
\end{equation*}

The family $\Gh(\alpha)$ contains a certain family of minimal elements, which we now define.
For $\kappa \in [h, \alpha]$ and $\lambda \in [0,1]$, we define the function
\begin{equation*}
    h_{\kappa, \lambda}(\theta) =
    \begin{cases}
        \kappa(1-\theta) &: \theta \geq \lambda\\
        \kappa(1-\lambda) &: \theta \leq \lambda.
    \end{cases}
\end{equation*}
This is the family of functions from \cite[Corollary~C]{zbmath:07937992}, and has the following minimality properties.
\begin{lemma}\label{l:omega-min}
    Let $0\leq h\leq\alpha$, $\kappa \in [0, \alpha]$, and $\lambda \in [0,1]$.
    Then the following hold:
    \begin{enumerate}[nl,r]
        \item $h_{\kappa,\lambda}\in\overline{\mathcal{G}}_0(\alpha)$.
        \item\label{i:2a} If $\gamma\in\mathcal{G}(\alpha)$ and $\omega = \Omega_h(\gamma)$ is such that $h_{\kappa,\lambda}(\lambda) \leq \omega(\lambda)$, then
            \begin{equation*}
                h_{\kappa,\lambda}(\theta) \leq \omega(\theta)\quad\text{for all}\quad 0\leq \theta \leq 1.
            \end{equation*}
        \item\label{i:3a} If $\psi \in \Gh(\alpha)$ is such that $h_{\kappa,\lambda}(\lambda) \leq \psi(\lambda)$, then
            \begin{equation*}
                h_{\kappa,\lambda}(\theta) \leq \psi(\theta)\quad\text{for all}\quad 0\leq \theta \leq 1.
            \end{equation*}
    \end{enumerate}
\end{lemma}
\begin{proof}
    It is immediate to see that $h_{\kappa,\lambda}\in\overline{\mathcal{G}}_0(\alpha)$.
    Also, \cref{i:2a} and \cref{i:3a} are trivial if $\lambda = 1$, since $h_{\kappa,1} = 0$ for all $\kappa$.
    Thus we may assume that $\lambda < 1$.

    Now we prove \cref{i:2a}.
    Let $\gamma\in\mathcal{G}(\alpha)$ and $\omega = \Omega_h(\gamma)$ be such that $h_{\kappa,\lambda}(\lambda) \leq \omega(\lambda)$.
    If $\omega(\lambda) = (1-\lambda)h$, then $\kappa = h$ and $h_{\kappa,\lambda}(\theta) \leq (1-\theta) h \leq \omega(\theta)$ for all $\theta$.

    Otherwise, get $\lambda' \leq \lambda$ such that
    \begin{equation*}
        \frac{\omega(\lambda)}{1-\lambda} = \frac{\gamma(\lambda')}{1-\lambda'}.
    \end{equation*}
    If $\theta \geq \lambda'$, since $h_{\kappa,\lambda}(\theta) \leq \kappa(1-\theta)$,
    \begin{equation*}
        \frac{\omega(\theta)}{1-\theta} \geq \frac{\gamma(\lambda')}{1-\lambda'} = \frac{\omega(\lambda)}{1-\lambda} \geq \frac{h_{\kappa,\lambda}(\lambda)}{1-\lambda} \geq \frac{h_{\kappa,\lambda}(\theta)}{1-\theta};
    \end{equation*}
    and if $\theta \leq \lambda'$, since $\omega \geq \gamma$ and $\gamma$ is decreasing,
    \begin{equation*}
        \omega(\theta) \geq \gamma(\theta) \geq \gamma(\lambda') = \left(\frac{1-\lambda'}{1-\lambda}\right)\omega(\lambda) \geq \omega(\lambda).
    \end{equation*}
    Since $h_{\kappa,\lambda}$ is constant on $[\theta,\lambda]$, it follows that $\omega(\theta) \leq h_{\kappa,\lambda}(\theta)$.
    This handles all cases of $\theta$, as required.

    Finally we prove \cref{i:3a}.
    Suppose $\psi \in \Gh(\alpha)$ is such that $h_{\kappa,\lambda}(\lambda) \leq \psi(\lambda)$.
    Since $\theta \mapsto \psi(\theta)/(1-\theta)$ is increasing and $\theta\mapsto h_{\kappa,\lambda}(\theta)/(1-\theta)$ is constant, $h_{\kappa,\lambda}(\theta) \geq \psi(\theta)$ for $\theta \geq \lambda$.
    Then since $\psi$ is decreasing and $h_{\kappa,\lambda}$ is constant on $[0,\lambda)$, $h_{\kappa,\lambda}(\theta) \geq \psi(\theta)$ for $\theta \leq \lambda$, as claimed.
\end{proof}
We now establish the following analogue of \cref{t:mono-sub-B}.
\begin{theorem}\label{t:mono-sub}
    Let $0 \leq h\leq \alpha$ and $\gamma\in\mathcal{G}(\alpha)$.
    Then $\Omega_h(\gamma)\in\Gh(\alpha)$, and
    \begin{equation*}
        \Omega_h(\gamma) = \inf\{\xi\in\Gh(\alpha): \xi \geq \gamma\}.
    \end{equation*}
    In particular, $\Omega_h\colon\mathcal{G}(\alpha) \to\Gh(\alpha)$ is a surjective idempotent map which is the identity on $\Gh(\alpha)$.
\end{theorem}
\begin{proof}
    For each $0 \leq \lambda \leq 1$, let $\kappa(\lambda)$ be chosen so that $h_{\kappa(\lambda),\lambda}(\lambda) = \gamma(\lambda)$.
    Certainly $\kappa(\lambda)\leq \alpha$, so $h_{\kappa,\lambda} \in \overline{\mathcal{G}}_0(\alpha)$ and therefore by the Arzelà--Ascoli theorem,
    \begin{equation*}
        \omega_0 \coloneqq \max_{0\leq \lambda \leq 1}h_{\kappa(\lambda),\lambda} \in \overline{\mathcal{G}}_0(\alpha).
    \end{equation*}
    Then writing
    \begin{equation*}
        \omega(\theta) = \max \{\omega_0(\theta), h(1-\theta)\},
    \end{equation*}
    it follows that $\omega \in \Gh(\alpha)$.
    Thus by \cref{l:omega-min}~\cref{i:3a},
    \begin{equation*}
        \omega = \inf\{\xi\in\Gh(\alpha): \xi \geq \gamma\};
    \end{equation*}
    and by \cref{l:omega-min}~\cref{i:2a}, $\omega_0 \leq \Omega_h(\gamma)$, and therefore $\omega \leq \Omega_h(\gamma)$.

    Therefore it remains to show that $\omega \geq \Omega_h(\gamma)$.
    Let $0\leq \theta < 1$ be arbitrary.
    If $\Omega_h(\gamma)(\theta) = h(1-\theta)$, clearly $\omega(\theta)\geq \Omega_h(\gamma)(\theta)$.
    Otherwise, by definition of $\Omega_h$, get $0 \leq \lambda \leq\theta$ such that
    \begin{equation*}
        \Omega_h(\gamma)(\theta) = \left(\frac{1-\theta}{1-\lambda}\right)\gamma(\lambda).
    \end{equation*}
    Then
    \begin{equation*}
        \Omega_h(\gamma)(\theta) = \left(\frac{1-\theta}{1-\lambda}\right) h_{\kappa(\lambda),\lambda}(\lambda) = h_{\kappa(\lambda),\lambda}(\theta) \leq \omega(\theta).
    \end{equation*}
    This completes the proof that $\omega = \Omega_h(\gamma)$, as required.
\end{proof}

\subsection{Equivalence of monotone subspace maps}
Finally, we show that the monotone subspaces $\Bh(\alpha)$ and $\Gh(\alpha)$ are equivalent.

Recall from \cref{p:Gamma-inv} that the function $\Gamma\colon\mathcal{B}(\alpha)\to\mathcal{G}(\alpha)$ is an order-preserving surjection with maximal right-inverse
\begin{equation*}
    \Gamma^{-1}(\gamma)(u,v) = u \gamma(v/u).
\end{equation*}
The heart of the matter is the observation that $\Gamma$ and $\Gamma^{-1}$ descend to maps on the monotone subspaces.
\begin{lemma}\label{l:gamma-descends}
    If $\psi\in\Bh(\alpha)$, then $\Gamma(\psi)\in\Gh(\alpha)$.
    Conversely, if $\gamma\in\Gh(\alpha)$, then $\Gamma^{-1}(\gamma) \in \Bh(\alpha)$.
\end{lemma}
\begin{proof}
    First, suppose $\psi\in\Bh(\alpha)$ is arbitrary: we show that $\gamma\coloneqq\Gamma(\psi)\in\Gh(\alpha)$.
    Since $\psi(u,0) \geq h u$, it is immediate that $\gamma(0) \geq h$.
    It remains to show that $\theta \mapsto \gamma(\theta)/(1-\theta)$ is increasing.
    Let $0\leq \lambda \leq \theta < 1$ and $u>0$ be arbitrary.
    Let $z\geq 0$ be chosen so that $\theta(u+z) = \lambda u + z$.
    We compute:
    \begin{equation*}
        \psi_{u+z}(\theta) = \frac{\psi(u+z,\lambda u +z)}{u+z} \geq \frac{\psi(u,\lambda u)}{u}\cdot\left(\frac{u}{u+z}\right) = \psi_u(\lambda)\cdot\left(\frac{1-\theta}{1-\lambda}\right).
    \end{equation*}
    Taking a limit supremum in $u$ proves that $\gamma\in\Gh(\alpha)$.

    Conversely, let $\gamma\in\Gh(\alpha)$ and let $\psi = \Gamma^{-1}(\gamma)$.
    We begin with the monotonicity property.
    Let $(u,v)\in\Delta$ and $z \geq 0$ be arbitrary.
    Since $\Psi(\gamma)$ is increasing,
    \begin{equation*}
        \psi(u+z, v+z) = (u+z)\cdot \gamma\left(\frac{v+z}{u+z}\right) \geq (u+z)\cdot\gamma\left(\frac{v}{u}\right)\cdot \left(\frac{1-\frac{v+z}{u+z}}{1-\frac{v}{u}}\right) = \psi(u,v).
    \end{equation*}
    Moreover, since $\gamma(0) \geq h$, if $0 \leq v \leq u$, then
    \begin{equation*}
        \psi(u,0) - \psi(v,0) = \gamma(0)(u-v) \geq h(u-v).
    \end{equation*}
    This proves that $\psi\in\Bh(\alpha)$, as required.
\end{proof}
Using \cref{l:gamma-descends} and the projection formulas in \cref{t:mono-sub-B} and \cref{t:mono-sub}, proving the equivalence of $\Bh(\alpha)$ and $\Gh(\alpha)$ is an algebraic formality.

\begin{theorem}\label{t:mono-sub-equiv}
    Let $0\leq h\leq \alpha$.
    Then $\Gamma\circ\Phi_h = \Omega_h \circ\Gamma$ as maps from $\mathcal{B}(\alpha)$ to $\Gh(\alpha)$.
    In particular, $\Gamma\colon\Bh(\alpha) \to\Gh(\alpha)$ is surjective.
\end{theorem}
\begin{proof}
    Let $\psi\in\mathcal{B}(\alpha)$ be arbitrary and write
    \begin{align*}
        A&\coloneqq\{\Gamma(g): g\in\Bh(\alpha)\text{ and }\psi \leq g\}\\
        B&\coloneqq \{f: f\in\Gh(\alpha)\text{ and }\Gamma(\psi) \leq f\}.
    \end{align*}
    By \cref{t:mono-sub-B} and \cref{t:mono-sub}, it suffices to prove that $A=B$.

    If $g\in\Bh(\alpha)$ and $\psi\leq g$, since $\Gamma$ is order-preserving, $\Gamma(\psi) \leq \Gamma(g)$.
    Since $\Gamma(\psi)\in\Gh(\alpha)$ by \cref{l:gamma-descends}, it follows that $A\subseteq B$.

    Conversely, suppose $f\in\Gh(\alpha)$ and $\Gamma(\psi) \leq f$.
    Write
    \begin{equation*}
        g = \max\{\overline{\Omega}_h(\psi), \Gamma^{-1}(f)\}.
    \end{equation*}
    Observe that $g\in\Bh(\alpha)$ by \cref{l:gamma-descends}, and of course $\psi\leq g$.
    Moreover, by maximality of $\Gamma^{-1}(f)$ from \cref{p:Gamma-inv}, $\Gamma(g) = f$.
    Therefore $B\subseteq A$, and equality holds.

    Finally, we recall from \cref{p:Gamma-inv} that $\Gamma\colon\mathcal{B}(\alpha)\to\mathcal{G}(\alpha)$ is surjective,
    Since $\Omega_h$ is surjective, it follows that $\Gamma\colon\Bh(\alpha) \to\Gh(\alpha)$ is surjective.
\end{proof}

\section{Dimensions of inhomogeneous attractors}\label{s:inhomog-2}
\subsection{Inhomogeneous attractors and separation conditions}\label{s:inhom-assump}
Fix a non-empty compact quasi-doubling metric space $X$ and a finite non-empty family of maps $\{f_i\}_{i\in\mathcal{I}}$ where each $f_i\colon X\to X$ is a strict Lipschitz contraction.
Then for a compact $F\subset X$, by the contraction mapping principle, there exists a unique non-empty compact $\Lambda_F \subset X$ satisfying the invariance relationship
\begin{equation*}
    \Lambda_F = F \cup \bigcup_{i\in\mathcal{I}}f_i(\Lambda_F).
\end{equation*}
We refer to the set $\Lambda_F$ as the \emph{inhomogeneous attractor} associated with the IFS $\{f_i\}_{i\in\mathcal{I}}$ and compact set $F$.

As is common in the study of iterated function systems, we require more assumptions on the maps $f_i$ in order to say something meaningful.
Let $\mathcal{I}^*$ denote the set of all finite words on $\mathcal{I}$.
For $\mtt{i} = (i_1,\ldots,i_n)$, we write $f_{\mtt{i}} = f_{i_1}\circ\cdots\circ f_{i_n}$.

We begin with regularity conditions.
First, let us introduce some terminology for contracting maps which do not distort space too much.
Let $(X,d)$ be a metric space and let $z\in\R$.
We define a new metric space $X_z$ with the same underlying set $X$ and with metric $d_z(x,y) = 2^{-z}d(x,y)$.
An easy computation shows that if $\diam X \leq 1$ and $z \geq 0$,
\begin{equation*}
    \beta_{X_z} = T_z(\beta_X)\quad\text{where}\quad
    T_z(\psi)(u,v) = \begin{cases}
        \psi(u-z,v-z) &: 0\leq z \leq v \leq u,\\
        \psi(u-z, 0) &: 0\leq v \leq z \leq u,\\
        0 &: 0\leq v\leq u \leq z.
    \end{cases}
\end{equation*}
Moreover, it is easy to check that $T_z$ maps $\mathcal{B}$ into itself, and similarly maps $\mathcal{B}(\alpha)$ into itself for all $\alpha \geq 0$.
In other words, at the level of branching functions, rescaling by a factor of $2^{-z}$ corresponds to the map $T_z$ given above.
\begin{definition}
    Let $X$ and $Y$ be metric spaces, let $\eta\in\mathcal{E}$, and let $z \geq 0$.
    Write $\eta_z(u) = \eta(u+z)$.
    We say that $f\colon X\to Y$ is \emph{$(\eta, z)$-contracting} if $f$ is a Lipschitz contraction and $\id\circ f\colon X \to Y_{-z}$ is $\eta_z$-uniform.
\end{definition}
This condition is somewhat weaker than requiring that $f$ is a bi-Lipschitz contraction with contraction ratio $2^{-z}$.
\begin{lemma}\label{l:ctr-equiv}
    Suppose $\diam X \leq 1$ and $f\colon X \to Y$ is surjective and Lipschitz contracting.
    Then $f$ is $(\eta, z)$-contracting for some $z \geq 0$ and $\eta\in\mathcal{E}$ if and only if
    \begin{equation*}
        |T_z(\beta_X)(u,v) - \beta_Y(u,v)| \leq \eta(u).
    \end{equation*}
\end{lemma}
\begin{proof}
    Since $f$ is Lipschitz contracting, $\diam Y \leq 1$ as well.
    Therefore, $f$ is $(\eta, z)$-contracting if and only if for all $(u,v)\in\Delta$
    \begin{equation*}
        |\beta_X(u,v) - \beta_{Y_{-z}}(u,v)| \leq \eta(u + z).
    \end{equation*}
    if and only if for all $(u,v)\in\Delta$
    \begin{equation*}
        |\beta_{X_z}(u,v) - \beta_Y(u,v)| \leq \eta(u).
    \end{equation*}
    Recalling that $\beta_{X_z} = T_z(\beta_X)$ since $z \geq 0$ and $\diam X \leq 1$ gives the claim.
\end{proof}

We can now state our main regularity condition.
\begin{definition}\label{d:md}
    We say that the IFS $\{f_i\}_{i\in\mathcal{I}}$ is \emph{minimally distorting} if there is a function $\rho\colon\mathcal{I}^*\to[0,\infty)$ and a function $\eta\in\mathcal{E}$ such that:
    \begin{enumerate}[nl,r]
        \item $\rho(\varnothing) = 0$ and $\rho(\mtt{i})>0$ for $\mtt{i} \neq \varnothing$.
        \item\label{i:subadd} There is a constant $A \geq 0$ such that for all $\mtt{i}, \mtt{j}\in\mathcal{I}$,
            \begin{equation*}
                0 \leq \rho(\mtt{i}\mtt{j}) - (\rho(\mtt{i}) + \rho(\mtt{j})) \leq A.
            \end{equation*}
        \item For all $\mtt{i}\in\mathcal{I}^*$, $f_{\mtt{i}}\colon F\to f_{\mtt{i}}(F)$ is $(\eta, \rho(\mtt{i}))$-contracting.
    \end{enumerate}
\end{definition}
Next, for $u>0$ and $x\in X$, we write
\begin{align*}
    \mathcal{I}^*(u) &= \{\mtt{i}\in\mathcal{I}^*:\rho(\mtt{i}^-) < u \leq \rho(\mtt{i})\}.\\
    \mathcal{F}(x, u) &= \{\mtt{i}\in\mathcal{I}^*(u):f_{\mtt{i}}(X)\cap B(x,2^{-u}) \neq \varnothing\}.
\end{align*}
Here, $\mtt{i}^-$ denotes the prefix of $\mtt{i}$ with length $|\mtt{i}| - 1$.
Note that every sufficiently long word has a exactly one prefix in $\mathcal{I}^*(u)$.

Now, let $\{f_i\}_{i\in\mathcal{I}}$ be minimally distorting with respect to a function $\rho$.
We define the \emph{critical exponent}
\begin{equation*}
    h = \lim_{u\to\infty}\frac{\log\#\mathcal{I}^*(u)}{u}.
\end{equation*}
The existence of the limit follows by a subadditivity argument using \cref{d:md}~\cref{i:subadd}.

We now state a non-concentration condition.
\begin{definition}\label{d:asymp}
    We say that the IFS $\{f_i\}_{i\in\mathcal{I}}$ is \emph{asymptotically bounded} if
    \begin{equation*}
        \lim_{u\to \infty}\frac{\log \sup_{x\in X}\#\mathcal{F}(x,u)}{u} = 0.
    \end{equation*}
\end{definition}

For the remainder of \cref{s:inhomog-2}, we make the following standing assumptions and fix some accompanying notation.
\begin{enumerate}
    \item We assume that $\{f_i\}_{i\in\mathcal{I}}$ is minimally distorting with functions $\rho$ and $\eta$, and asymptotically bounded with critical exponent $h$.
    \item We also fix a compact $F\subset X$, and consider the attractor $\Lambda = \Lambda_F$.
    \item Each map $f_i$ fixes some $x_i\in X$, so we may assume without loss of generality that $x_i\in F$ for all $i$; this not change the attractor $\Lambda$ and only changes $\beta_F$ by a fixed constant factor.
        A convenient consequence is that $F \cap f_{\mtt{i}}(F)\neq \varnothing$ for all $\mtt{i}\in\mathcal{I}^*$.
    \item We assume that $\diam X \leq 1$.
\end{enumerate}

\subsection{Geometric lemmas}
We begin by observing that at resolution $r$, the attractor $\Lambda$ essentially looks like a union of copies of $F$ arranged in a tree structure governed by the underlying IFS.
\begin{lemma}\label{l:r-nbhd}
    For all $u \geq 0$, $\Lambda$ is a subset of the closed $r$-neighbourhood of
    \begin{equation*}
        \bigcup_{\substack{\mtt{i}\in\mathcal{I}^*\\\rho(\mtt{i}) > r}}f_{\mtt{i}}(F).
    \end{equation*}
\end{lemma}
\begin{proof}
    It is easy to check that
    \begin{equation*}
        \Lambda = \overline{\bigcup_{\mtt{i}\in\mathcal{I}^*}f_{\mtt{i}}(F)}.
    \end{equation*}
    Moreover, if $\rho(\mtt{i}) \leq r$, since $f_{\mtt{i}}(F) \cap F \neq \varnothing$, $f_{\mtt{i}}(F)$ is contained in the $(r\diam X)$-neighbourhood of $F$.
    Since we assumed that $\diam X \leq 1$, we are done.
\end{proof}
We next note (using asymptotic boundedness) that we can always find large well-separated subsets of sets of cylinders in $\mathcal{I}^*(u)$.
\begin{lemma}\label{l:large-sub}
    Let $u > 0$ and $\mathcal{W} \subset \mathcal{I}^*(u)$ be non-empty.
    Then there exists a subset $\mathcal{W}' \subset \mathcal{W}$ such that $\log\#\mathcal{W}' = \log\#\mathcal{W} + o(u)$ and moreover for each $x\in X$,
    \begin{equation}\label{e:disj}
        \#\{\mtt{i}\in\mathcal{W}':f_{\mtt{i}}(F) \cap B(x, 2^{-u}) \neq \varnothing\} \leq 1.
    \end{equation}
\end{lemma}
\begin{proof}
    We first note that $\diam f_{\mtt{i}}(F) \leq 2^{-\rho(\mtt{i})} \leq 2^{-u}$.
    In particular, $f_{\mtt{i}}(F) \subset B(x_{\mtt{i}}, 2^{-u})$ for any fixed choice $x_{\mtt{i}}\in f_{\mtt{i}}(F)$.
    Moreover, by asymptotic boundedness, $B(x_{\mtt{i}}, 2\cdot 2^{-u})$ intersects $2^{o(u)}$ cylinders in $\mathcal{W}$.

    Now, construct $\mathcal{W}'$ inductively.
    Begin with some $\mtt{i}\in\mathcal{W}$, add $\mtt{i}$ to $\mathcal{W}'$, and delete all $\mtt{j}\in\mathcal{W}$ such that $f_{\mtt{j}}(F) \cap B(x_{\mtt{i}}, 2\cdot 2^{-u}) \neq \varnothing$.
    This removes at most $o(u)$ cylinders from $\mathcal{W}$.
    Repeat for new $\mtt{i}\in\mathcal{W}$ until $\mathcal{W}$ is empty.

    Let us check that $\mathcal{W}'$ satisfies the desired properties.
    Clearly, $2^{o(u)}\cdot\#\mathcal{W}' \geq \#\mathcal{W}$.
    Moreover, if $x\in X$ and $f_{\mtt{i}}(F) \cap B(x, 2^{-u})\neq\varnothing$, then $B(x, 2^{-u}) \subset B(x_{\mtt{i}}, 2\cdot 2^{-u})$ and therefore \cref{e:disj} holds for the ball $B(x, 2^{-u})$, as required.
\end{proof}

\subsection{Cylinder counting}
We write $\rho_{\min} = \min_{i\in\mathcal{I}}\rho(i)$ and $\rho_{\max} = \max_{i\in\mathcal{I}}\rho(i)$.
By assumption, $0 < \rho_{\min} \leq \rho_{\max}< \infty$.
Note that if $\mtt{i}\in\mathcal{I}^*$ and $\mtt{j}\in\mathcal{I}^*$, then we have the crude bound
\begin{equation*}
    |\mtt{j}|\rho_{\min} \leq \rho(\mtt{i} j) - \rho(\mtt{i}) \leq A + |\mtt{j}|\rho_{\max}.
\end{equation*}

We begin with discretization lemma for $\mathcal{I}^*$.
\begin{lemma}\label{l:disc}
    We have
    \begin{equation*}
        \mathcal{I}^* = \bigcup_{k=0}^\infty\mathcal{I}^*(k\rho_{\min}).
    \end{equation*}
\end{lemma}
\begin{proof}
    Let $\mtt{i}\in\mathcal{I}^*$ be arbitrary: then $\mtt{i}\in\mathcal{I}^*(2^{-u})$ for any $u$ such that $\rho(\mtt{i}^-) < 2^{-u} \leq \rho(\mtt{i})$.
    Since $\rho(\mtt{i}) \geq \rho(\mtt{i}^-) + \rho_{\min}$, we can choose $u = k\rho_{\min}$ for some integer $k\geq 0$.
\end{proof}

Now, we establish our core cylinder count for the number of cylinders which can intersect a given ball, even when the cylinders and balls potentially have very different size.
More precisely, given $v\geq 0$, $z \geq 0$, and $x\in\Lambda$, set
\begin{equation*}
    \mathcal{J}(v, z, x) = \left\{\mtt{i}\in\mathcal{I}^*(z):f_{\mtt{i}}(F)\cap B(x, 2^{-v})\neq\varnothing\right\}.
\end{equation*}
We now have the following key lemma.
\begin{lemma}\label{l:cyl-count}
    Let $v\geq 0$ and $z \geq 0$.
    Then
    \begin{equation*}
        \log\sup_{x\in \Lambda}\#\mathcal{J}(v,z,x) =
        \begin{cases}
            o(v) &: z \leq v\\
            h(z-v) + o(v) + o(z-v) &: v\leq z.
        \end{cases}
    \end{equation*}
\end{lemma}
\begin{proof}
    If $0\leq z \leq v$, then for all $x\in\Lambda$, we have $\mathcal{J}(v,z,x) \subset \mathcal{F}(v, x)$ and the claim follows by asymptotic boundedness.

    Otherwise, let $v \leq z$.
    We begin with the upper bound.
    Fix $x\in\Lambda$ and $\mtt{j}\in\mathcal{F}(v, x)$, and consider the set
    \begin{equation*}
        \mathcal{P}(\mtt{j})\coloneqq\{\mtt{k}\in\mathcal{I}^*: \mtt{j}\mtt{k} \in \mathcal{I}^*(z)\}.
    \end{equation*}
    Unpacking the definitions and applying \cref{d:md}~\cref{i:subadd},
    \begin{equation*}
        \rho(\mtt{k}^-) < z - v \leq \rho(\mtt{k}) + A + \rho_{\max}.
    \end{equation*}
    In particular, there is a number $\ell$ depending only on the global parameters such that
    \begin{equation*}
        \mathcal{P}(\mtt{j}) \subset \bigcup_{k=0}^\ell\mathcal{I}^*(z-v+k\rho_{\min}).
    \end{equation*}
    Finally, if $\mtt{i}\in\mathcal{J}(v,z,x)$, then $\mtt{i}$ has a unique prefix $\mtt{j}\in\mathcal{F}(z, x)$, in which case $\mtt{i}\in\mathcal{P}(\mtt{j})$.
    Therefore
    \begin{align*}
        \log\#\mathcal{J}(v,z,x) &\leq \log\sum_{\mtt{j}\in \mathcal{F}(z, x)} \#\mathcal{P}(\mtt{j})\\
                                 &\leq \max_{k=0,\ldots,\ell}\log\#\mathcal{I}^*(z-v+k\rho_{\min}) + o(v)\\
                                 &\leq h(z-v) + o(v) + o(z-v).
    \end{align*}
    In the second inequality, we used asymptotic boundedness to bound $\#\mathcal{F}(z,x)$.

    We now conclude with the lower bound.
    If $h=0$ there is nothing to prove.
    Otherwise, let $0<\varepsilon<h$ be arbitrary and by definition of $h$ get $u_0 >0$ so that for all $u \geq u_0$,
    \begin{equation*}
        \log\#\mathcal{I}^*(u) \geq u(h-\varepsilon).
    \end{equation*}
    For $z-v\leq u_0$, there is nothing to prove.
    Otherwise, let $x\in\Lambda$ be arbitrary and get $\mtt{j}\in\mathcal{I}(v)$ such that $x\in f_{\mtt{j}}(X)$.
    Since $\diam X \leq 1$, $\mtt{j}\in \mathcal{J}(v,v,x)$.
    Now if $\mtt{k}\in\mathcal{I}^*(z-v)$, then
    \begin{align*}
        \rho(\mtt{j}\mtt{k}) &\geq \rho(\mtt{j}) + \rho(\mtt{k}) \geq z\\
        \rho(\mtt{j}\mtt{k}^-) & \leq A + \rho(\mtt{j}) + \rho(\mtt{k}^-) \leq A + \rho_{\max} + z.
    \end{align*}
    Therefore, $\mtt{j}\mtt{k}$ has a unique prefix $\mtt{i}\in\mathcal{I}^*(z)$, and since $f_{\mtt{j}}(X)\subset B(x, 2^{-v})$, $\mtt{i}\in\mathcal{J}(v,z,x)$.
    Moreover, by the second equation in the above display combined with \cref{l:disc}, there is some number $\ell$ depending only on the global parameters such that each such prefix $\mtt{i}$ appears from at most $\ell$ words $\mtt{k}$.
    Therefore
    \begin{equation*}
        \log\#\mathcal{J}(v,z,x) \geq \log\#\mathcal{I}^*(z-v) - \log\ell \geq (z-v)(h-\varepsilon) -\log\ell.
    \end{equation*}
    Since $\varepsilon>0$ was arbitrary, the lower bound follows.
\end{proof}

\subsection{Proof of the dimension formula}

In this section, we prove our main result concerning inhomogeneous attractors.
\begin{theorem}\label{t:in-dim}
    Let $\{f_i\}_{i\in\mathcal{I}}$ be asymptotically bounded and minimally distorting, with critical exponent $h$, and let $F\subset X$ be non-empty and compact.
    Then $\beta_{\Lambda} = \Phi_h(\beta_F)$.
\end{theorem}
\begin{proof}
    Let $(u,v)\in\Delta$ be arbitrary.
    By \cref{l:r-nbhd}, decomposing $[0, u]$ into $O(u)$ intervals of width $\rho_{\min}$,
    \begin{equation*}
        \beta_\Lambda(u,v) = \sup_{0 \leq z \leq u}\omega(z) + O(\log u)
    \end{equation*}
    where
    \begin{equation*}
        \omega(z)\coloneqq\sup_{x\in\Lambda}\log N_{2^{-u}}\left(B(x, 2^{-v})\cap \bigcup_{\mtt{i}\in\mathcal{J}(v,z,x)} f_{\mtt{i}}(F)\right).
    \end{equation*}
    For the remainder of the proof, we estimate $\omega(z)$.
    We consider two cases depending on the value of $z$.

    First, suppose $0\leq z \leq v$.
    Since $\log\#\mathcal{J}(v,z,x) = o(z)$ by \cref{l:cyl-count},
    \begin{align*}
        \omega(z) &= \sup_{x\in\Lambda}\max_{\mtt{i}\in\mathcal{J}(v,z,x)}\log N_{2^{-u}}\left(B(x, 2^{-v})\cap f_{\mtt{i}}(F)\right) + o(u)\\
                  &= \sup_{x\in\Lambda}\max_{\mtt{i}\in\mathcal{J}(v,z,x)}\beta_{f_{\mtt{i}}(F)}(u,v) + o(u)\\
                  &= \beta_{F}(u-z,v-z) + o(u).
    \end{align*}
    In the last line, we used the fact that $f_{\mtt{i}}$ is $(\eta, \rho(\mtt{i}))$-contracting, recalling that $z - \rho_{\max} \leq \rho(\mtt{i}) \leq z$.

    Otherwise, $v \leq z \leq u$.
    Again, $\sup_{x\in \Lambda}\log\#\mathcal{J}(v,z,x) = h(z-v) + o(u)$ by \cref{l:cyl-count}, and therefore
    \begin{align*}
        \omega(z) &\leq \log \#\mathcal{J}(v,z,x) + \max_{\mtt{i}\in\mathcal{J}(v,z,x)}\beta_{f_{\mtt{i}}(F)}(u, v) + o(u)\\
                  &= h(z-v) + \beta_F(u-z, 0) + o(u).
    \end{align*}
    To get the lower bound, apply \cref{l:large-sub} to get $\mathcal{J}'\subset \mathcal{J}(v,z,x)$ with $\log\#\mathcal{J}' = \log\#\mathcal{J}(v,z,x) + o(u)$ such that any ball $B(x, 2^{-u})$ intersects at most $1$ set $f_{\mtt{i}}(F)$ for $\mtt{i}\in\mathcal{J}'$.
    But then
    \begin{equation*}
        \omega(z) \geq \log \#\mathcal{J}' + \beta_F(u-z, 0) + o(u) \geq h(z-v) + \beta_F(u-z,0) + o(u).
    \end{equation*}
    Recalling the definition of $\Phi_h$, this completes the proof, as required.
\end{proof}

\begin{acknowledgements}
    VO is supported by National Research, Development and Innovation Office - NKFIH, Project K142169 NKFI KKP144059 “Fractal geometry and applications” and by the Hungarian Research Network through the HUN-REN BME Stochastics Research group.
    AR is supported by the Research Council of Finland via Tuomas Orponen's project \emph{Approximate incidence geometry}, grant no.\ 355453.
\end{acknowledgements}
\appendix
\section{Alternative definition of \texorpdfstring{$\Gh(\alpha)$}{Gh(α)}}\label{s:alt}
In this section, we show that the subadditivity condition in the definition of $\Gh(\alpha)$ is redundant.
\begin{definition}
    For $0 \leq h \leq \alpha$, we let $\Gh(\alpha)$ denote the functions $\gamma\colon[0,1]\to [0,\alpha]$ such that:
    \begin{enumerate}[nl,r]
        \item $\gamma(0) \geq h$ and $\gamma(1) = 0$,
        \item $\gamma$ is decreasing and $\alpha$-Lipschitz, and
        \item $\theta\mapsto\gamma(\theta)/(1-\theta)$ is increasing.
    \end{enumerate}
\end{definition}
The following lemma was observed in \cite[Lemma~2.6]{zbmath:07937992}, but we give the short details here for completeness.
\begin{lemma}
    Let $\alpha \geq 0$.
    Then $\Gh(\alpha) \subset \mathcal{G}(\alpha)$.
\end{lemma}
\begin{proof}
    Let $\gamma\in\Gh(\alpha)$.
    It suffices to check \cref{d:G-def}~\cref{i:G-subadd}: that is, for $\lambda,\theta \in [0,1]$ that
    \begin{equation*}
        \gamma(\lambda\theta) \leq \gamma(\theta) + \theta \gamma(\lambda).
    \end{equation*}
    If $\theta \leq \lambda$, using the increasing property \cref{i:incr} with $\lambda\theta \leq \lambda$ and then $\theta \leq \lambda$,
    \begin{equation*}
        \gamma(\lambda\theta) - \gamma(\theta)\leq\frac{1-\lambda\theta}{1-\theta}\cdot \gamma(\theta) - \gamma(\theta) = \theta\cdot \frac{1-\lambda}{1-\theta}\gamma(\theta) \leq \theta \gamma(\lambda),
    \end{equation*}
    and if $\lambda \leq \theta$, using \cref{i:incr} with $\lambda \theta \leq \lambda$ and $\lambda \leq \theta$,
    \begin{equation*}
        \gamma(\lambda\theta) - \gamma(\theta) \leq \frac{1-\lambda \theta}{1-\lambda} \gamma(\lambda) - \frac{1-\theta}{1-\lambda}\gamma(\lambda) = \theta \gamma(\lambda).
    \end{equation*}
    It follows that $\gamma\in\mathcal{G}(\alpha)$.
\end{proof}

\section{Conformal IFSs are minimally distorting}\label{s:conf}
In this section, we show that conformal IFSs are minimally distorting.

Following \cite{zbl:0852.28005}, let $X$ be a compact connected subset of $\R^d$ with the Euclidean norm and let $\mathcal{I}$ be a non-empty finite index set.
Fix a family of injective strictly contracting maps $S_i \colon X \to X$ for $i \in \mathcal{I}$.

\begin{definition}\label{d:cifs}
    We say that the IFS $\{S_i\}_{i\in\mathcal{I}}$ is \emph{conformal} if the following additional properties are satisfied:
    \begin{enumerate}[r]
        \item\label{i:conformal} \emph{Conformality}:
            There exists an open, bounded, connected subset $V \subset \R^d$ such that $X \subset V$ and such that for each $i \in \mathcal{I}$, $S_i$ extends to a conformal $C^{1+\varepsilon}$ diffeomorphism on $V$.
        \item\label{i:bdp} \emph{Bounded distortion}:
            There exists $K\geq 1$ such that $\norm{S_{\mtt{i}}'(x)} \leq K\snorm{S_{\mtt{i}}'(y)}$ for all $x,y \in V$ and $\mtt{i} \in \mathcal{I}^*$.
            Here, $S_{\mtt{i}}'(x)$ denotes the Jacobian of the map $S_{\mtt{i}}$ at $x$ and $\norm{\cdot}$ denotes the spectral matrix norm.
    \end{enumerate}
\end{definition}
Any conformal IFS satisfies the minimally distorting property.
\begin{proposition}
    If $\{S_i\}_{i\in\mathcal{I}}$ is conformal, then it is minimally distorting.
\end{proposition}
\begin{proof}
    For $\mtt{i}\in\mathcal{I}^*$, define
    \begin{equation*}
        \rho(\mtt{i}) = -\log\sup_{x\in X}\norm{S_{\mtt{i}}'(x)}\qquad\text{for}\qquad\mtt{i}\in\mathcal{I}^*.
    \end{equation*}
    Since $S_{\varnothing}$ is the identity map, $\rho(\varnothing) = 0$; since the IFS is strictly contracting, $\rho(\mtt{i}) > 0$ for $\mtt{i}\neq\varnothing$.
    Then by the chain rule, sub-multiplicativity of the matrix norm, and the bounded distortion property with constant $K$, it follows that
    \begin{equation*}
        0 \leq \rho(\mtt{i}\mtt{j}) - (\rho(\mtt{i}) + \rho(\mtt{j}))\leq \log K.
    \end{equation*}
    Finally, it is known that there is a constant $C>0$ depending only on the IFS such that
    \begin{equation*}
        C^{-1} \leq \frac{\norm{S_{\mtt{i}}(x) - S_{\mtt{i}}(y)}}{\rho(\mtt{i})\cdot\norm{x-y}} \leq C;
    \end{equation*}
    see for example the discussion in \cite[pp.~111--112]{zbl:0852.28005} or \cite[Lemma~2.9]{zbl:1547.28014}.
    Applying \cref{p:pres}, this shows that $S_{\mtt{i}}$ is $(O(1), \rho(\mtt{i}))$-contracting, which completes the proof of the minimal distortion property.
\end{proof}

\end{document}